\documentclass[10pt]{amsart}

\usepackage{amsmath,amsthm}
\usepackage{amssymb}
\usepackage{euscript}

\input epsf

\numberwithin{equation}{section}

\newtheorem{thm}[equation]{Theorem}
\newtheorem{prop}[equation]{Proposition}
\newtheorem{cor}[equation]{Corollary}
\newtheorem{lemma}[equation]{Lemma}

\theoremstyle{definition}
\newtheorem{defn}[equation]{Definition}

\newcommand{\cn}{\mathcal{C}}
\newcommand{\bcn}{\overline{\mathcal{C}}}
\newcommand{\hm}{\mathcal{H}}
\newcommand{\bhm}{\overline{\mathcal{H}}}
\newcommand{\nhm}{\textsl{H}}
\newcommand{\nbhm}{\overline{\textsl{H}}}

\DeclareMathOperator{\spn}{span}
\DeclareMathOperator{\sgn}{sign}
\DeclareMathOperator{\rnk}{rank}
\DeclareMathOperator{\Ker}{Ker}

\DeclareMathOperator{\Img}{Im}

\newcommand{\my}{\Phi}
\newcommand{\myd}{(\my + d)}
\newcommand{\mydpr}{(\my' + d')}

\newcommand{\1}{\mathbf{1}}
\newcommand{\x}{\mathbf{x}}
\newcommand{\fa}{\mathbf{a}}
\newcommand{\fb}{\mathbf{b}}
\newcommand{\A}{\mathbf{A}}
\newcommand{\M}{\mathbf{M}}
\newcommand{\indI}{\mathcal{I}}
\newcommand{\indJ}{\mathcal{J}}
\newcommand{\isom}{\rho}

\begin{document}

\title[Endomorphism of Khovanov Invariant]{An endomorphism of the Khovanov 
invariant}
\author{Eun Soo Lee}

\begin{abstract}
We construct an endomorphism of the Khovanov invariant to prove H-thinness 
and pairing phenomena of the invariants for alternating links. As a 
consequence, it follows that the Khovanov invariant of an oriented nonsplit 
alternating link is determined by its Jones polynomial, signature, and the 
linking numbers of its components.
\end{abstract}

\date{\today}
\subjclass{57M27}
\keywords{Khovanov invariant, H-thin, alternating links}

\address{Department of Mathematics, MIT, Cambridge, MA 02139, U.S.A.}
\email{eslee@alum.mit.edu}

\maketitle

\bigskip

\section{Introduction}
\label{scn:intro}

Khovanov invariant is a cohomology theory for oriented links with values in 
graded abelian groups, and specializes to the Jones polynomial by taking 
graded Euler characteristic of those cohomology groups 
(theorem \ref{thm:Jones}). Khovanov \cite{K} constructed the invariant in a 
search of connections between combinatorial invariants and differential 
geometric invariants of 3 and 4 dimensional manifolds. He interpreted his 
coboundary map as the image of a functor from the category of two 
dimensional cobordisms between one dimensional manifolds to the category 
of $\mathbb{Z}[c]$-modules.

The Khovanov invariant specialized by setting $c=0$ and tensoring with 
$\mathbb{Q}$ (will be just called the Khovanov invariant from now on) has 
been computed by Bar-Natan \cite{B}\cite{B2} for the prime knots with up to 
11 crossings. From Bar-Natan's data, two conjectures \cite{B}\cite{G} on 
the values of Khovanov invariant for alternating knots were formulated by 
Bar-Natan, Garoufalidis, and Khovanov.
The conjectures (theorems \ref{thm:conj2} and \ref{thm:conj1}) imply that 
the Khovanov invariant of an alternating knot determines and is determined 
by its Jones polynomial and signature.

The following is the theorem in \cite{K} which states that the Khovanov 
invariant specializes to the Jones polynomial. The Khovanov invariant of a 
(relatively) oriented link $L$ in rational coefficients is denoted by 
$\hm(L)$ following \cite{K}, and is defined in section \ref{scn:inv}. 
Its associated polynomial is denoted by $Kh(L)$ as it is in \cite{B}.
\[ Kh(L)(t,q) \buildrel \textrm{def} \over = \sum t^i q^j \dim \hm^{i,j}(L) \]

\begin{thm}[\cite{K}]
\label{thm:Jones}
For an oriented link $L$, the graded Euler characteristic
\[ \sum_{i,j \in \mathbb{Z}} (-1)^i q^j \dim \hm^{i,j}(L) \]
of the Khovanov invariant $\hm(L)$ of $L$ is equal to $(q^{-1} + q)$ times 
the Jones polynomial $V(L)$ of $L$.
\[ \sum_{i,j \in \mathbb{Z}} (-1)^i q^j \dim \hm^{i,j}(L) = (q^{-1} + q) 
V(L)_{\sqrt{t} = -q} \]

In terms of the associated polynomial $Kh(L)$,
\[ Kh(L)(-1, q) = (q^{-1} + q) V(L)_{\sqrt{t} = -q} \textrm{ .} \]
\end{thm}

The following two theorems are the conjectures in \cite{B} proved in this 
paper.

\begin{thm}[Conjecture 2 in \cite{B} and \cite{G}]
\label{thm:conj2}
For any alternating knot $L$, the Khovanov invariants $\hm^{i,j}(L)$ of $L$ 
are supported in two lines
\[ j = 2 i - \sigma(L) \pm 1 \textrm{ .} \]

In other words, the equality
\[ Kh(L)(t,q) = q^{-\sigma(L)}(q^{-1} \cdot A(tq^2) + q \cdot B(tq^2)) \]
holds for some polynomials $A$ and $B$, where $\sigma(L)$ is the signature 
of $L$.
\end{thm}

Theorem \ref{thm:conj2}, in fact, holds for any (relatively) oriented 
nonsplit alternating link $L$. (See theorem \ref{thm:conj2-ext}.)

\begin{defn}[H-thinness \cite{K2}]
A diagram/knot/link $D$ is \emph{H-thin} if its Khovanov invariant $\hm(D)$ 
(or $\bhm(D)$ if $D$ is a diagram) is supported in two diagonal lines as in 
theorem \ref{thm:conj2} up to a shift of the grading.
\end{defn}

Theorem \ref{thm:conj2} implies that any nonsplit alternating link is H-thin.

\begin{thm}[Conjecture 1 in \cite{B} and \cite{G}]
\label{thm:conj1}
For an alternating knot $L$, its Khovanov invariants $\hm^{i,j}(L)$ of degree 
difference $(1,4)$ are paired except in the 0th cohomology group.

More precisely, in terms of the polynomial $Kh(L)$, the equality
\[ Kh(L)(t,q) = q^{-s} (q^{-1} + q) + (q^{-1} + t q^2 \cdot q) \cdot C(t,q) \]
holds for some integer $s$ and some polynomial $C$.
\end{thm}

Theorem \ref{thm:conj1} was extended to (relatively) oriented nonsplit 
alternating links in theorem \ref{thm:conj1-ext}.

Combining theorems \ref{thm:conj2} and \ref{thm:conj1}, we can write $Kh(L)$ as
\[ Kh(L)(t,q) = q^{-\sigma(L)} \{ (q^{-1} + q) + (q^{-1} + t q^2 \cdot q) 
\cdot Kh'(L)(t q^2)\} \]
for some polynomial $Kh'(L)$.

As it is discussed in \cite{B} and \cite{G}, theorems \ref{thm:conj2} and 
\ref{thm:conj1} with theorem \ref{thm:Jones} imply that the Khovanov invariant, 
or equivalently the associated polynomial $Kh(L)$, of an alternating knot $L$ 
is completely determined by the Jones polynomial and the signature of $L$.

Fortunately, that is not the case for nonalternating knots. A counterexample 
can be found in \cite{B2}: $10_{136}$ and $11^n_{92}$ both have signature $-2$ 
and the same Jones polynomial, but their Khovanov invariants do not agree.

The organization of the coming sections is as follows. Section
\ref{scn:inv} consists of a brief summary of the Khovanov invariant. Section 
\ref{scn:2nd} is devoted to our
proof of theorem \ref{thm:conj2}. In section \ref{scn:1st}, an 
endomorphism of Khovanov invariant is defined and used to prove
\ref{thm:conj1}.

We follow \cite{L} and \cite{Br} for basic notions in knot theory and graph 
theory, and \cite{K} for notations and terminologies related to the Khovanov 
invariant. We only need a relative orientation to define the Khovanov 
invariant, so an orientation and oriented can be read as a relative 
orientation and relatively oriented.

\bigskip

\section{Khovanov Invariant}
\label{scn:inv}

In this section, the construction and some properties of Khovanov
invariant in \cite{K} are summarized. Khovanov's original
construction is more general, but we will concentrate on a
specialized case with coefficients in $\mathbb{Q}$. The interested
reader should read \cite{K}.

\subsection{Construction}

\subsubsection{Cubes of diagrams}

Let $L$ be an oriented link and $D$ be its diagram, a regular projection 
of $L$ together with the information of relative height at each double point.
A double point of $D$ can be resolved in two ways.

\[ \epsfxsize=2.92in\epsfbox{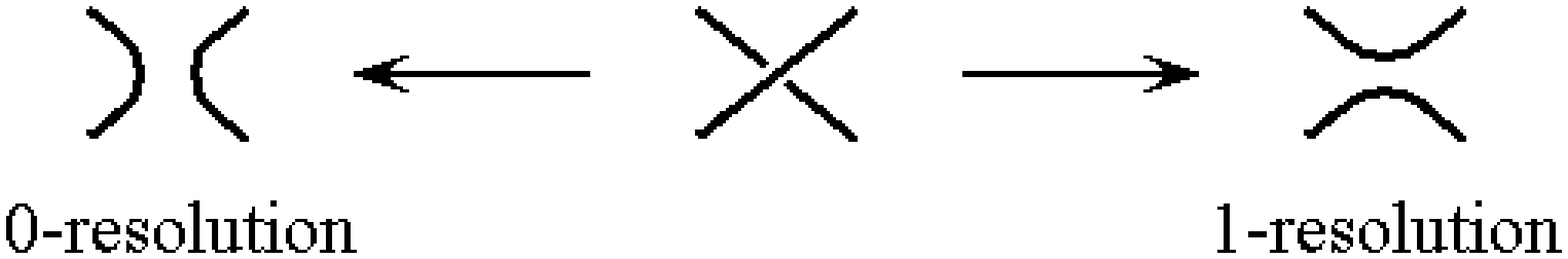} \]

Let $\indI$ be the set of double points of $D$. Each subset $\indJ$ of 
$\indI$ corresponds to a complete resolution $D(\indJ)$ of $D$ in which 
points in $\indJ$ are resolved to their 1-resolutions, points not in 
$\indJ$ are to their 0-resolutions.
Regard those subsets of $\indI$ as vertices.
For each pair of vertices $\indJ$ and $\indJ'$ satisfying
$\indJ \subset \indJ'$ and $|\indJ' - \indJ| = 1$, there is a directed 
edge from $\indJ$ to $\indJ'$. A directed cube is constructed.

\subsubsection{Cubes of modules}

Let $\A = \mathbb{Q} \1 \oplus \mathbb{Q} \x$ be a two dimensional module 
over $\mathbb{Q}$ with a multiplication $m$, a comultiplication $\Delta$, 
a unit $\iota$, and a counit $\epsilon$ defined as
\begin{eqnarray*}
m(\1 \otimes \1) &=& \1 \\
m(\1 \otimes \x) = m(\x \otimes \1) &=& \x \\
m(\x \otimes \x) &=& 0 \\
\Delta(\1) &=& \1 \otimes \x + \x \otimes \1 \\
\Delta(\x) &=& \x \otimes \x \\
\iota (1) &=& \1 \\
\epsilon (\1) &=& 0 \\
\epsilon (\x) &=& 1 \textrm{ .}
\end{eqnarray*}

For each vertex $\indJ$ of a cube, assign a tensor product of as many copies 
of $\A$ as the number of components of $D(\indJ)$, and denote it by 
$\M_{\indJ}(D)$. There is a one-to-one correspondence between those copies 
of $\A$ and the components of $D(\indJ)$.

\subsubsection{Chain complexes}

A chain complex can be constructed from the cube of modules. Its $i$-th chain 
group is a direct sum of all the modules over vertices of $i$ elements.
\[ \bcn^i(D) = \bigoplus_{|\indJ| = i} \M_{\indJ}(D) \]

To define the coboundary map $d$, choose an ordering of $\indI$ - the set of 
crossings of $D$, and regard $\indJ \subset \indI$ as an ordered 
$|\indJ|$-tuple of its elements in the chosen order instead of just a subset 
of $\indI$.

For a homogeneous element $x \in \M_{\indJ}(D)$, $dx$ lies in the sum of all 
the modules over those vertices which are end-points of the directed edges 
from $\indJ$.
\[ dx \in \bigoplus_{\indJ \subset \indJ', |\indJ' - \indJ| = 1} 
\M_{\indJ'}(D) \hspace{.5in} \textrm{ for } x \in \M_{\indJ}(D) \]
Each homogeneous component of $d$ is defined in the way that 
$m: \A \otimes \A \rightarrow \A$ is applied to 
corresponding modules if two components merge into one, and 
$\Delta: \A \rightarrow \A \otimes \A$ is if one component splits to two.
If the ordered $(|\indJ| + 1)$-tuple $\indJ$ followed by the element in 
$\indJ' - \indJ$ is an odd permutation of the ordered 
$(|\indJ| + 1)$-tuple $\indJ'$, $-m$ or $-\Delta$, instead of $m$ or 
$\Delta$, is used for the $\M_{\indJ}(D) \rightarrow \M_{\indJ'}(D)$ 
component of $d$. With this choice of signs, $d$ satisfies $d^2 = 0$.

\subsubsection{Relation to TQFT}
\label{sbsbscn:cat}

The algebra $\A$ above is a Frobenius algebra and it is related to a two 
dimensional topological quantum field theory. There is a functor $F$ that 
maps one dimensional manifolds consisting of $n$ disjoint simple closed 
curves to $\A^{\otimes n}$'s, and cobordisms in the following figure to
\[ F(S^1_2) = m, F(S^2_1) = \Delta, F(S^1_0) =\iota, F(S^0_1) = \epsilon, 
F(S^2_2) = (\textrm{permutation}), F(S^1_1) = id \textrm{ .} \]

\[ \epsfxsize=3.23in\epsfbox{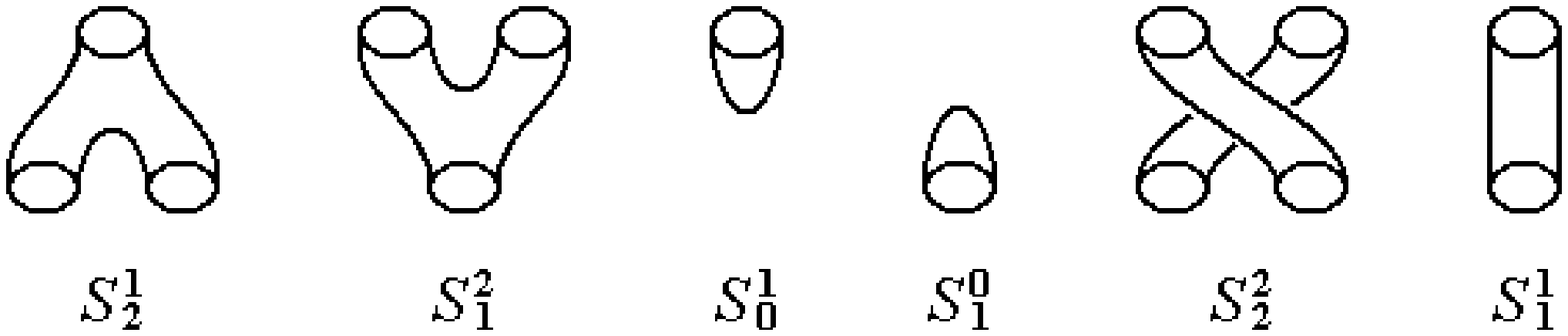} \]

The unit, counit, (co)associativity, (co)commutativity, together with the 
identity
\[ \Delta \circ m = (m \otimes id) \circ (id \otimes \Delta) \]
ensures well-definedness of $F$. (See \cite{K}.) 

\[ \epsfxsize=3.73in\epsfbox{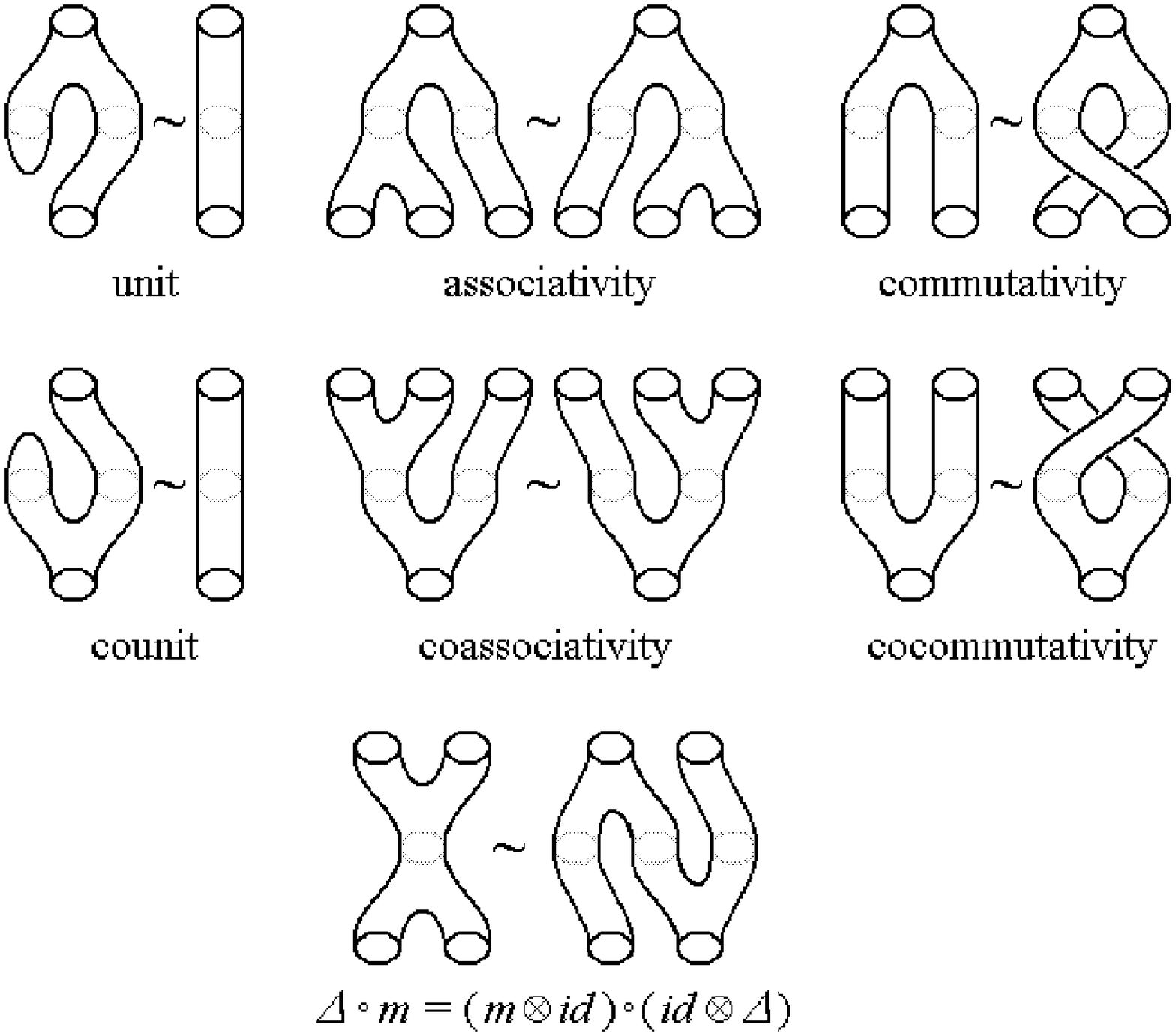} \]

From the viewpoint of the previous sections, well-definedness of
$F$ implies $d^2=0$. The following figure tabulates all the
possible relative locations of two crossings and the associated
surfaces obtained by continuous change of resolutions of the two
crossings one after the other.

\[ \epsfxsize=3.95in\epsfbox{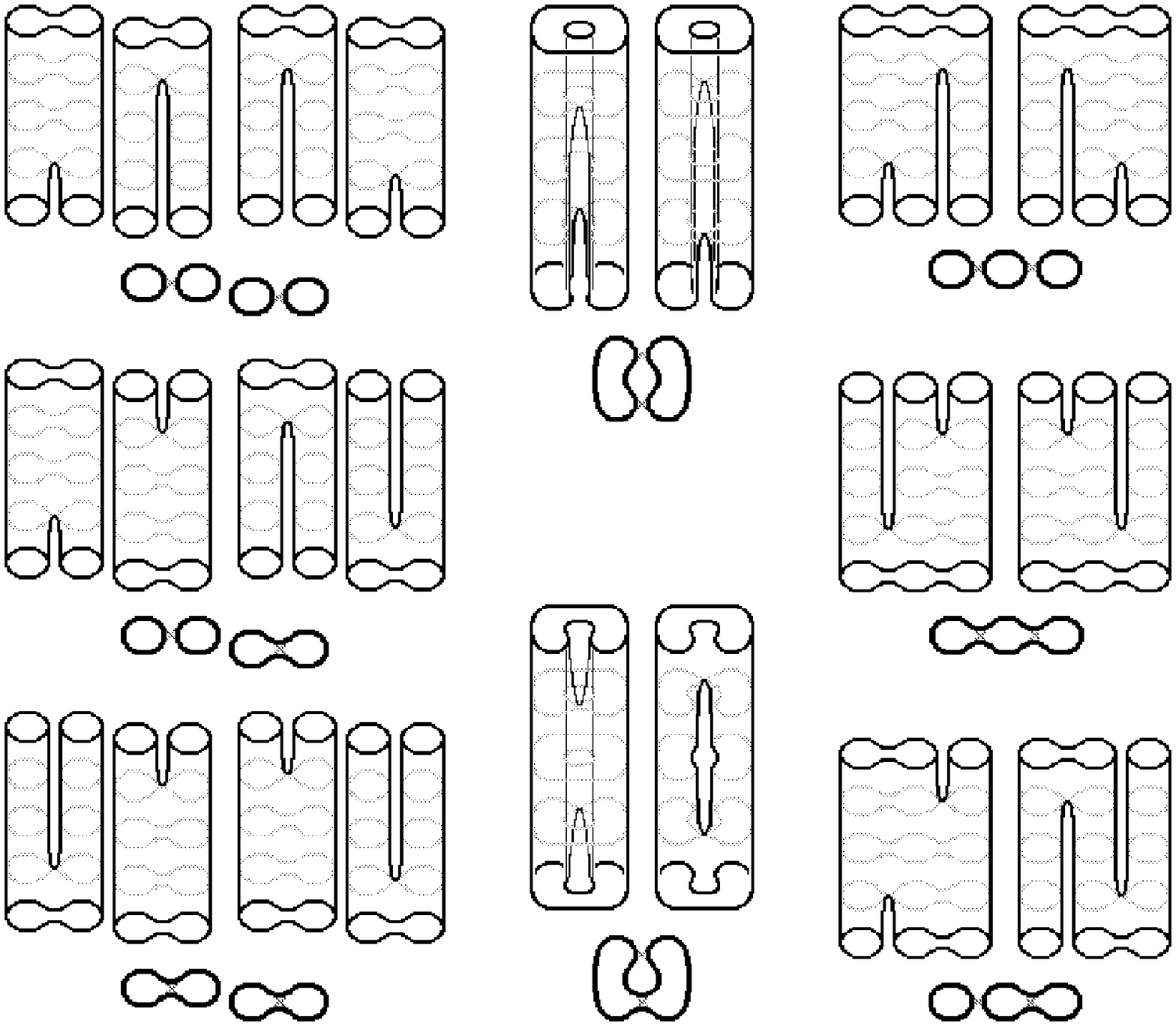} \]

For $d^2$ to be equal to 0, those in the first two columns from the left only 
require $m$ and $\Delta$ to be (co)commutative, the top right one requires 
$m$ to be associative, the middle right one requires $\Delta$ to be 
coassociative, and the bottom right one requires the identity
$\Delta \circ m = (m \otimes id) \circ (id \otimes \Delta)$.
The unit and counit can be dropped if we are concerned only about $d^2=0$.
(See \cite{A} for more discussion of Frobenius algebra associated to a 
two-dimensional TQFT.)

\subsubsection{Grading on the chain groups}

$\A$ is a graded module. $\x$ is of degree $-1$, $\1$ is of degree 1. 
$\M_{\indJ}(D)$ has a grading induced from that of $\A$. Note that $m$ 
and $\Delta$ are maps of degree $-1$ with respect to these gradings.

The chain group $\bcn^i(D)$ above has a grading shifted from those of 
$\M_{\indJ}(D)$:
\[ \bcn^i(D) = \bigoplus_{|\indJ|=i} \M_{\indJ}(D)\{-i\}\textrm{ .} \]
$\M\{k\}$ means a module identical to $\M$ with a shifted grading. 
An element of degree $j$ in $\M$ is of degree $j-k$ in $\M\{k\}$.

Define $\bcn^{i,j}(D)$ as the degree $j$ component of $\bcn^{i}(D)$. 
Due to the shifts, $d$ is now degree preserving, so that $\bhm^{i}(D)$ 
is also decomposed as $\bigoplus_j \bhm^{i,j}(D)$.

A chain complex $\cn(D)$ is defined from $\bcn(D)$ with the orientation 
of $D$ taken into account.

For each crossing of $D$, a sign is given as below (Note that this is 
opposite to the sign convention in \cite{K}.)

\[ \epsfxsize=1.91in\epsfbox{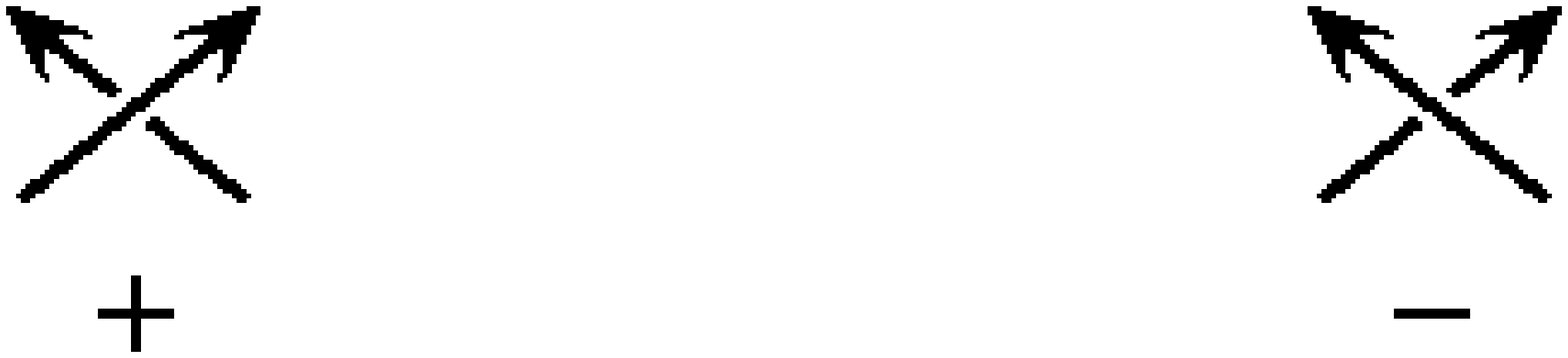} \]

Let $x(D)$ and $y(D)$ be the number of negative crossings and positive 
crossings of $D$, respectively. $\cn(D)$ is defined as
\[ \cn(D) = \bcn(D)[x(D)]\{2x(D)-y(D)\} \]
with the same coboundary map $d$. Square brackets indicate a shift of 
the indices of chain groups, i.e.,
\[ \cn^{i,j}(D) = \bcn^{i+x(D), j+2x(D)-y(D)}(D)\textrm{ .} \]

\subsubsection{Example}
\label{sbsbscn:trefoil}

Here is an illustration of what had happened in the previous
sections to the following diagram $D$ of the lefthanded trefoil.

\[ \epsfxsize=4.22in\epsfbox{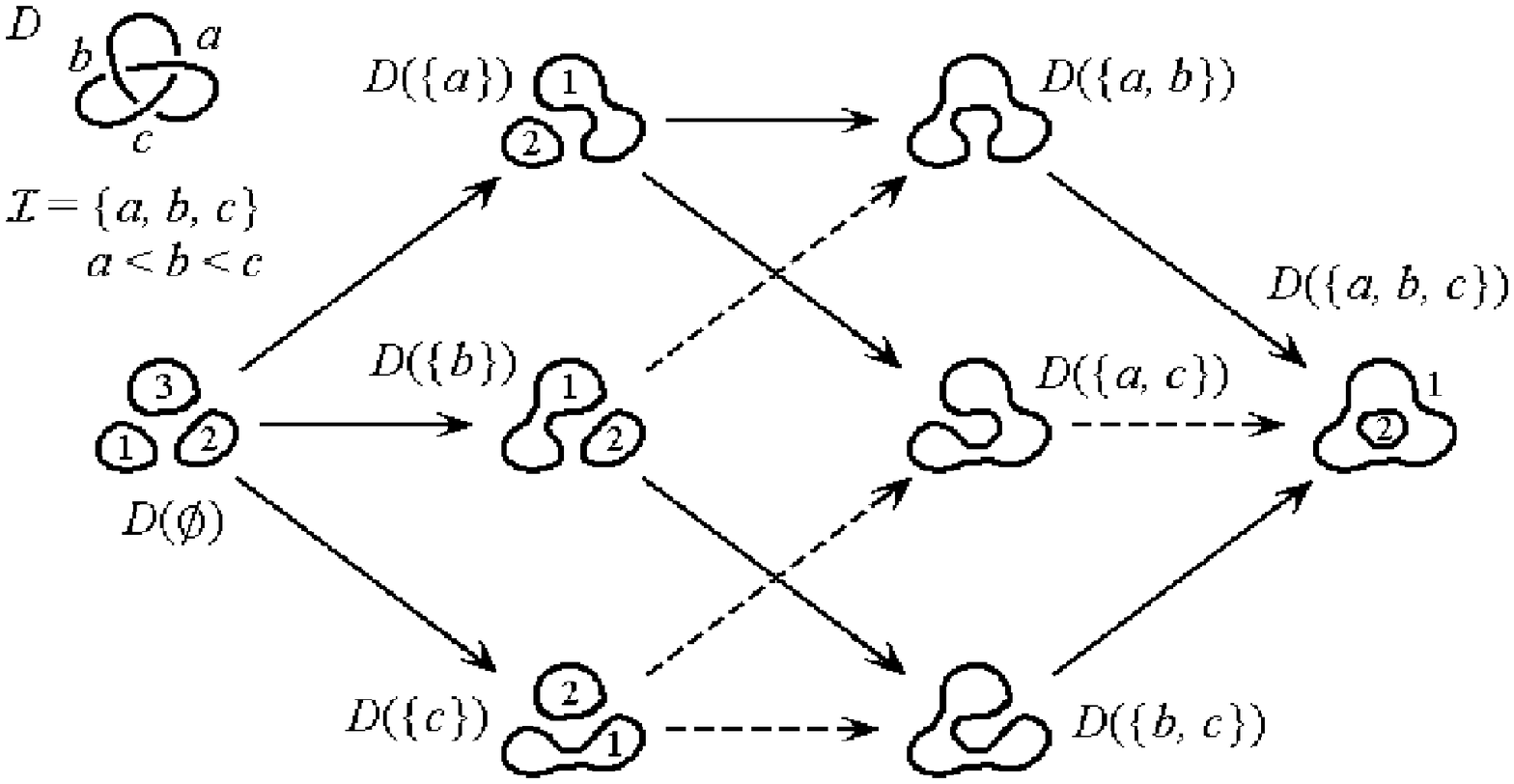} \]

Numbers for components of a resolution indicate which piece of $\A$ 
corresponds to which component, dotted edges indicate the places where 
$-m$ or $-\Delta$ should be used.

\begin{eqnarray*}
\bcn^0(D) &=& \M_{\emptyset}(D) = \A \otimes \A \otimes \A \\
\bcn^1(D) &=& \M_{\{a\}}(D) \oplus \M_{\{b\}}(D) \oplus \M_{\{c\}}(D) 
= (\A \otimes \A) \oplus (\A \otimes \A) \oplus (\A \otimes \A) \\
\bcn^2(D) &=& \M_{\{a,b\}}(D) \oplus \M_{\{a,c\}}(D) \oplus \M_{\{b,c\}}(D) 
= \A \oplus \A \oplus \A \\
\bcn^3(D) &=& \M_{\{a,b,c\}}(D) = \A \otimes \A
\end{eqnarray*}

\begin{eqnarray*}
\x \otimes \x \otimes \x \in \bcn^{0, -3}(D) & \buildrel d \over \longmapsto 
& \begin{bmatrix}0 \\ 0\\ 0\end{bmatrix} \in \bcn^{1, -3}(D) \\
\1 \otimes \x \otimes \1 \in \bcn^{0, 1}(D) & \buildrel d \over \longmapsto 
& \begin{bmatrix}\x \otimes \1 \\ \1 \otimes \x \\ \x \otimes \1 \end{bmatrix} 
\in \bcn^{1, 1}(D) \\
\begin{bmatrix}\1 \otimes \x \\ \x \otimes \1 \\ 0 \end{bmatrix} 
\in \bcn^{1, 1}(D) & \buildrel d \over \longmapsto 
& \begin{bmatrix}0 \\ \x \\ \x \end{bmatrix} \in \bcn^{2, 1}(D)\\
\begin{bmatrix}0 \\ \1 \\ 0 \end{bmatrix} \in \bcn^{2, 3}(D) 
& \buildrel d \over \longmapsto & - \1 \otimes \x - \x \otimes \1 
\in \bcn^{3, 3}(D)
\end{eqnarray*}

\begin{eqnarray*}
\bhm^{0, -3}(D) &=& \frac{\spn \{ \x \otimes \x \otimes \x \} }{ \{0\} }\\
\bhm^{1,1}(D) &=& \frac{\spn \{ 
\begin{bmatrix}\1 \otimes \x \\ \1 \otimes \x \\ \1 \otimes \x \end{bmatrix}, 
\begin{bmatrix}\1 \otimes \x \\ \x \otimes \1 \\ \x \otimes \1 \end{bmatrix}, 
\begin{bmatrix}\x \otimes \1 \\ \1 \otimes \x \\ \x \otimes \1 \end{bmatrix}, 
\begin{bmatrix}\x \otimes \1 \\ \x \otimes \1 \\ \1 \otimes \x \end{bmatrix} 
\} }{ \spn \{ 
\begin{bmatrix}\1 \otimes \x \\ \x \otimes \1 \\ \x \otimes \1 \end{bmatrix}, 
\begin{bmatrix}\x \otimes \1 \\ \1 \otimes \x \\ \x \otimes \1 \end{bmatrix}, 
\begin{bmatrix}\x \otimes \1 \\ \x \otimes \1 \\ \1 \otimes \x \end{bmatrix} 
\} }\\
\bhm^{3,3}(D) &=& \frac{ \spn \{ \1 \otimes \x, \x \otimes \1 \} }{ 
\spn \{ \1 \otimes \x + \x \otimes \1 \} }\\
\bhm^{3,5}(D) &=& \frac{ \spn \{ \1 \otimes \1 \} }{ \{0\} } \\
\bhm^{i, j}(D) &=&
\begin{cases}
\mathbb{Q} & \textrm{ for } (i, j) = (0, -3), (1, 1), (3, 3), 
\textrm{ or } (3, 5) \\
0 & \textrm{ otherwise }
\end{cases}
\end{eqnarray*}

$x(D) =3 $ and $y(D) =0$ for this diagram $D$. (A knot or a knot diagram has 
only one relative orientation.) Hence,
\[ \hm^{i,j}(D) = \bhm^{i+3, j+6}(D) =
\begin{cases}
\mathbb{Q} & \textrm{ for } (i, j) = (-3, -9), (-2, -5), (0, -3), 
\textrm{ or } (0, -1) \\
0 & \textrm{ otherwise }
\end{cases} \]

\subsection{Invariance}

To define $\hm(L)$ as $\hm(D)$, we need to see invariance of
$\hm(D)$ under a change of ordering of $\indI$ and under the
Reidemeister moves. We will just state isomorphisms. Detailed
proofs can be found in chapter 5 of \cite{K}.

\subsubsection{Change of ordering}

If the $|\indJ|$-tuple $\indJ$ with respect to one ordering of $\indI$ is an 
even permutation of $\indJ$ with respect to the other ordering, $\M_{\indJ}$ 
component of our isomorphism between them is the identity. Otherwise, it is 
minus identity.

\subsubsection{Reidemeister moves} 
\label{sbsbscn:move}

\[ \epsfxsize=4.02in\epsfbox{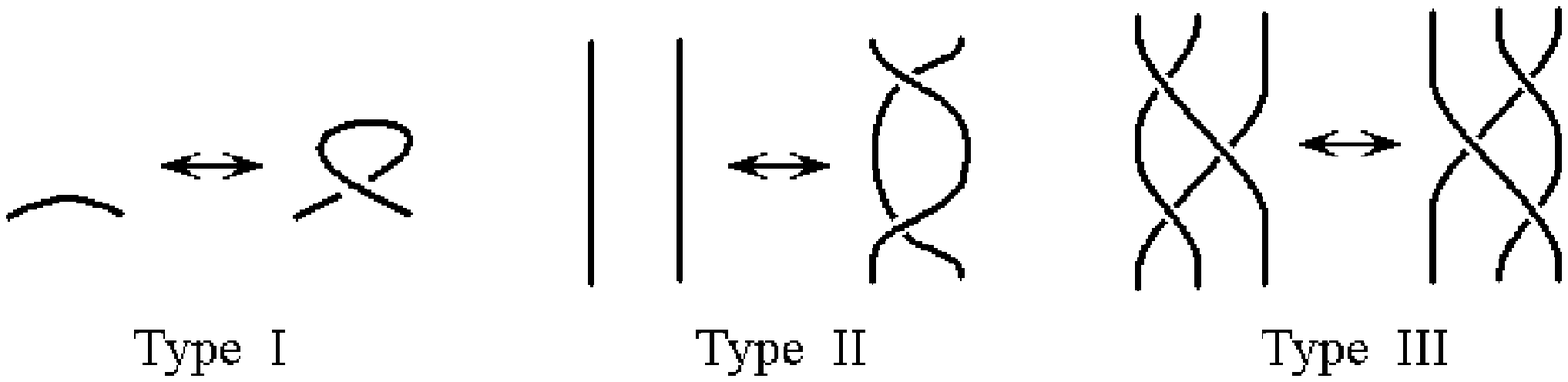} \]

[Type I]
\[ \epsfxsize=2.95in\epsfbox{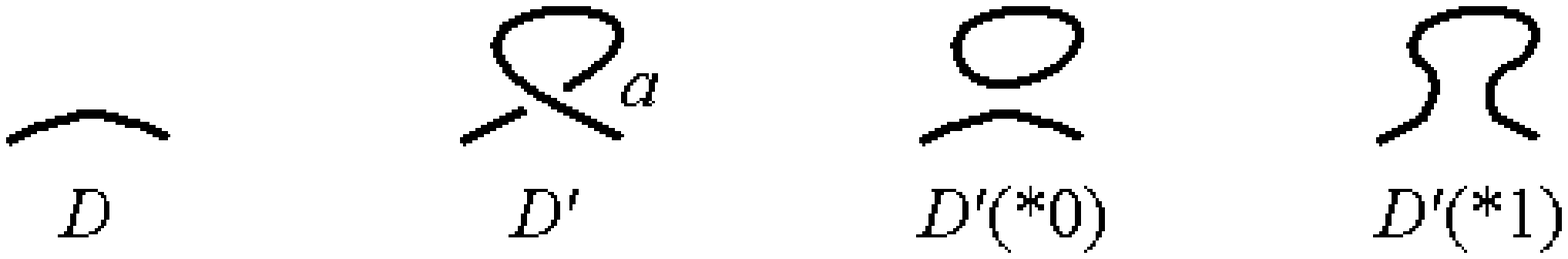} \]

Let $a$ be the crossing which appears only in $D'$. The set $\indI'$ of 
crossings of $D'$ is $\indI$, the set of crossings of $D$, followed by 
$a$ as an ordered $|\indI'|$-tuple. Let $D'(*0)$ and $D'(*1)$ denote $D'$ 
with only its last crossing (that is $a$) resolved to its 0- and 
1-resolutions, respectively.

As a group, $\bcn(D')$ is a direct sum of $\bcn(D'(*0))$ and 
$\bcn(D'(*1))[-1]\{-1\}$. Denote the part of the coboundary map $d'$ 
on $\bcn(D')$ that maps from $\bcn(D'(*0))$ to $\bcn(D'(*1))[-1]\{-1\}$ 
by $d'_{0 \rightarrow 1}$, and the coboundary maps on $\bcn(D(*0))$ 
and $\bcn(D(*1))$ by $d'_0$ and $d'_1$, so that
\[ d'(y + z) = d'_0 (y) + d'_{0 \rightarrow 1}(y) - d'_1(z) \]
for $y \in \bcn(D'(*0))$ and $z \in \bcn(D'(*1))[-1]\{-1\}$. Similar 
notation should be comprehended similarly from now on.

Define
\[ X_1 = \Ker d'_{0 \rightarrow 1} \]
and
\[ X_2 = \{ y \otimes \1 + z | y \in \bcn(D), z \in \bcn(D'(*1))[-1]\{-1\} \} 
\textrm{ .} \]
Here $\bcn(D'(*0))$ is identified with $\bcn(D) \otimes \A$.

$X_1$ and $X_2$ are subcomplexes of $\bcn(D')$, $\bcn(D')$ is decomposed as 
$X_1 \oplus X_2$ as a chain complex, $X_2$ is acyclic, and
\begin{eqnarray*}
\isom:& X_1 &\longrightarrow \; \bcn(D)\{1\} \\
&y \otimes \1 + z \otimes \x &\longmapsto \; z
\end{eqnarray*}
induces an isomorphism between $\hm(D')$ and $\hm(D)$.

[Type II]
\[ \epsfxsize=4.05in\epsfbox{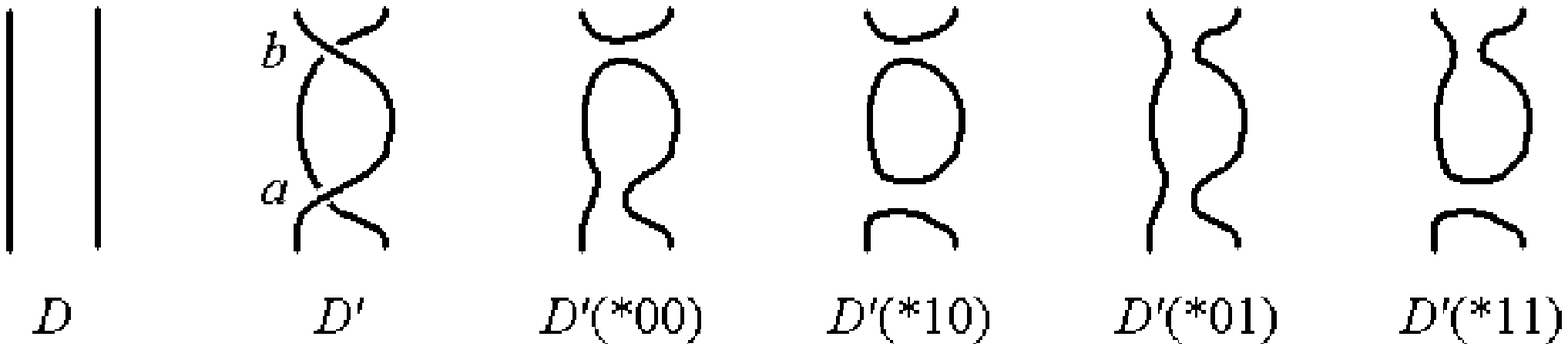} \]

As before, the set $\indI'$ of crossings of $D'$ is $\indI$, the set of 
crossings of $D$, followed by $a$, then $b$ as an ordered $|\indI'|$-tuple.

This time,
\begin{eqnarray*}
X_1 &=& \{ z + \alpha(z) | z \in \bcn(D'(*01))[-1]\{-1\} \} \\
X_2 &=& \{ z + d'y | z, y \in \bcn(D'(*00)) \} \\
X_3 &=& \{ z + y \otimes \1 | z, y \in \bcn(D'(*11))[-2]\{-2\}\}
\end{eqnarray*}
where $\alpha (z) = - d'_{01 \rightarrow 11}(z) \otimes \1 
\in \bcn(D'(*10))[-1]\{-1\} \approx \bcn(D'(*11))[-1]\{-1\} \otimes \A$.

Then, $\bcn(D')$ is a direct sum of its subcomplexes $X_1$, $X_2$, and $X_3$, 
$X_2$ and $X_3$ are acyclic, and
\begin{eqnarray*}
\isom: &\bcn^i(D)[-1]\{-1\} \approx \bcn^i(D'(*01))[-1]\{-1\}
& \longrightarrow \; X_1 \cap \bcn^i(D') \\
&z& \longmapsto \; (-1)^i(z + \alpha(z))
\end{eqnarray*}
induces an isomorphism between $\hm(D)$ and $\hm(D')$.

[Type III]
\[ \epsfxsize=4.73in\epsfbox{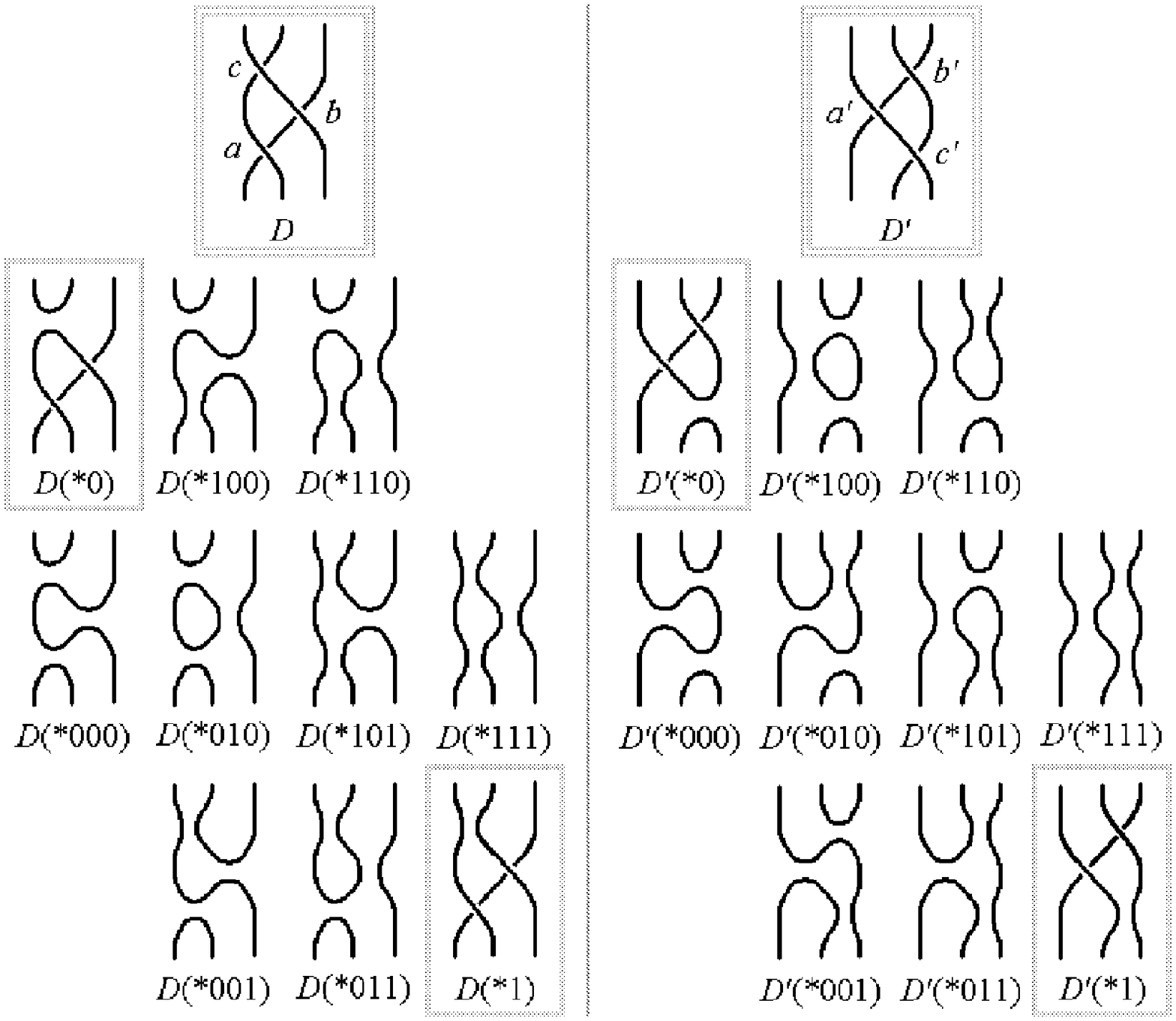} \]

Again, $a, b, c$ and $a', b', c'$ are the last three elements in $\indI$ 
and $\indI'$, and the others are in the same order.

Define $\alpha, \beta, \alpha', \beta'$ as
\begin{eqnarray*}
\alpha: &\bcn(D(*110))[-2]\{-2\}& \longrightarrow \; 
\bcn(D(*010))[-1]\{-1\} \approx \bcn(D(*110))[-1]\{-1\} \otimes \A \\
&z& \longmapsto \; z \otimes \1 \\
\beta: &\bcn(D(*100))[-1]\{-1\}& \longrightarrow \; 
\bcn(D(*010))[-1]\{-1\} \\
&z& \longmapsto \; \alpha d_{100 \rightarrow 110} (z) \\
\alpha': &\bcn(D'(*110))[-2]\{-2\}& \longrightarrow \; 
\bcn(D'(*100))[-1]\{-1\} \approx \bcn(D'(*110))[-1]\{-1\} \otimes \A \\
&z& \longmapsto \; z \otimes \1 \\
\beta': &\bcn(D'(*010))[-1]\{-1\}& \longrightarrow \; 
\bcn(D'(*100))[-1]\{-1\} \\
&z& \longmapsto \; - \alpha' d'_{010 \rightarrow 110} (z) \textrm{ .}
\end{eqnarray*}

$\bcn(D)$ and $\bcn(D')$ can be decomposed into their subcomplexes as below.
\begin{eqnarray*}
\bcn(D) &=& X_1 \oplus X_2 \oplus X_3 \\
X_1 &=& \{x + \beta(x) + y | x \in \bcn(D(*100))[-1]\{-1\}, 
y \in \bcn(D(*1))[-1]\{-1\} \} \\
X_2 &=& \{x + dy | x,y \in \bcn(D(*000))\} \\
X_3 &=& \{\alpha (x) + d \alpha (y)| x, y \in \bcn(D(*110))[-2]\{-2\} \} \\
\bcn(D') &=& Y_1 \oplus Y_2 \oplus Y_3 \\
Y_1 &=& \{x + \beta'(x) + y | x \in \bcn(D'(*010))[-1]\{-1\}, y \in 
\bcn(D'(*1))[-1]\{-1\} \} \\
Y_2 &=& \{x + d'y | x,y \in \bcn(D'(*000))\} \\
Y_3 &=& \{\alpha' (x) + d' \alpha' (y)| x, y \in \bcn(D'(*110))[-2]\{-2\} \}
\end{eqnarray*}

As before, $X_2, X_3, Y_2, Y_3$ are acyclic,
$\bcn(D(*100))[-1]\{-1\}$ and $\bcn(D'(*010))[-1]\{-1\}$, 
$\bcn(D(*1))[-1]\{-1\}$ and $\bcn(D'(*1))[-1]\{-1\}$ are naturally 
isomorphic, and $X_1$ is isomorphic to $Y_1$ via
\[ \isom: x + \beta(x) + y \longmapsto x + \beta'(x) + y \textrm{ .} \]

\subsection{Properties}

The following results are proved in \cite{K} and will be used in
the coming sections.

\begin{prop}
For an oriented $n$ component link diagram $D$,
\[ \cn^{i,j}(D) = 0 \]
unless $j \equiv n (\textrm{mod }2) $.
\end{prop}

\begin{cor}
For an oriented $n$ component link $L$,
\[ \hm^{i,j}(L) = 0 \]
unless $j \equiv n (\textrm{mod }2) $.
\end{cor}

\begin{prop}
For a disjoint union $D \sqcup D'$ of two oriented link diagrams $D$ and $D'$,
\[ \cn(D \sqcup D') = \cn(D) \otimes \cn(D') \textrm{ .} \]
\end{prop}

\begin{cor}
For a disjoint union $L \sqcup L'$ of two oriented links $L$ and $L'$,
\[ \hm(L \sqcup L') = \hm(L) \otimes \hm(L') \textrm{ .} \]
\end{cor}

\begin{prop}
The Frobenius algebra $(\A, m, \Delta, \iota, \epsilon)$ is isomorphic to its 
dual algebra $(\A^*, \Delta^*, m^*, \epsilon^*, \iota^*)$.
\end{prop}

\begin{prop}
Let $D^!$ be the mirror image of an oriented link diagram $D$. The complex 
$\cn(D^!)$ is isomorphic to the dual of $\cn(D)$.
\end{prop}

\begin{cor}
\label{cor:mirror}
For an oriented link $L$ and its mirror image $L^!$,
\[ \hm(L^!) \cong (\hm(L))^* \textrm{ .} \]
\end{cor}

\begin{thm}[Theorem \ref{thm:Jones}]
For an oriented link $L$, the graded Euler characteristic
\[ \sum_{i,j \in \mathbb{Z}} (-1)^i q^j \dim \hm^{i,j}(L) \]
of the Khovanov invariant $\hm(L)$ of $L$ is equal to $(q^{-1} + q)$ times 
the Jones polynomial $V(L)$ of $L$.
\[ \sum_{i,j \in \mathbb{Z}} (-1)^i q^j \dim \hm^{i,j}(L) = (q^{-1} + q) 
V(L)_{\sqrt{t} = -q} \]

In terms of the associated polynomial $Kh(L)$,
\[ Kh(L)(-1, q) = (q^{-1} + q) V(L)_{\sqrt{t} = -q} \textrm{ .} \]
\end{thm}

\bigskip

\section{H-thinness of Alternating Links}
\label{scn:2nd}

In this section, we prove theorem \ref{thm:conj2}. The proof is
based on induction on the number of crossings. We will show that
the support of the Khovanov invariant of a nonsplit oriented
alternating link is included in the union of the supports for two
such links with fewer crossings, then that the two lines of the two
supports agree.

\subsection{Exact sequences}

\begin{thm}
\label{thm:exact}
The chain complexes $\bcn(D)$, $\bcn(D(*0))$, and $\bcn(D(*1))[-1]\{-1\}$ 
form a short exact sequence
\[ 0 \rightarrow \bcn(D(*1))[-1]\{-1\} \rightarrow
\bcn(D) \rightarrow \bcn(D(*0)) \rightarrow 0
 \]
with degree preserving maps, so that $\bhm(D)$ is an extension
of the kernel and cokernel of the connecting map $\delta$ as a
bigraded $\mathbb{Q}$-module.
\[ \cdots \rightarrow \bhm^{i-1} (D(*0)) \buildrel \delta \over \rightarrow 
\bhm^{i-1}(D(*1))\{-1\} \rightarrow \bhm^i(D) \rightarrow \bhm^i (D(*0)) 
\buildrel \delta \over \rightarrow \bhm^i(D(*1))\{-1\} \rightarrow \cdots \]

In particular, the support of $\bhm(D)$ is included in the union of the 
support of $\bhm(D(*0))$ and $\bhm(D(*1))[-1]\{-1\}$.
\end{thm}

\begin{proof}
We've already seen that $\bcn(D)$ is decomposed as $\bcn(D(*0)) \oplus 
\bcn(D(*1))[-1]\{-1\}$, and the coboundary map $d$ can be written as
\[ d(y + z) = d_0 y + d_{0 \rightarrow 1} y - d_1 z \textrm{ .} \]

Now, it is easy to see that
\[ 0 \rightarrow \bcn(D(*1))[-1]\{-1\} \rightarrow
\bcn(D) \rightarrow \bcn(D(*0)) \rightarrow 0
 \]
is a short exact sequence of chain complexes (after a little
adjustment of sign), and that $\delta : \bhm^i(D(*0)) 
\to \bhm^{i+1}(D(*1))[-1]\{-1\}$ 
is induced by $d_{0 \rightarrow 1}$.
\end{proof}

\subsection{Properties of black and white coloring of alternating link 
diagrams}

Let $D$ be a link diagram.
For brevity of the statements to follow, let us think of diagrams on $S^2$ 
rather than on $\mathbb{R}^2$. The regions of $S^2$ divided by $D$ can be 
colored black and white in checkerboard fashion.

At each crossing, a coloring of nearby regions falls into one of the two 
following patterns.

\[ \epsfxsize=2.36in\epsfbox{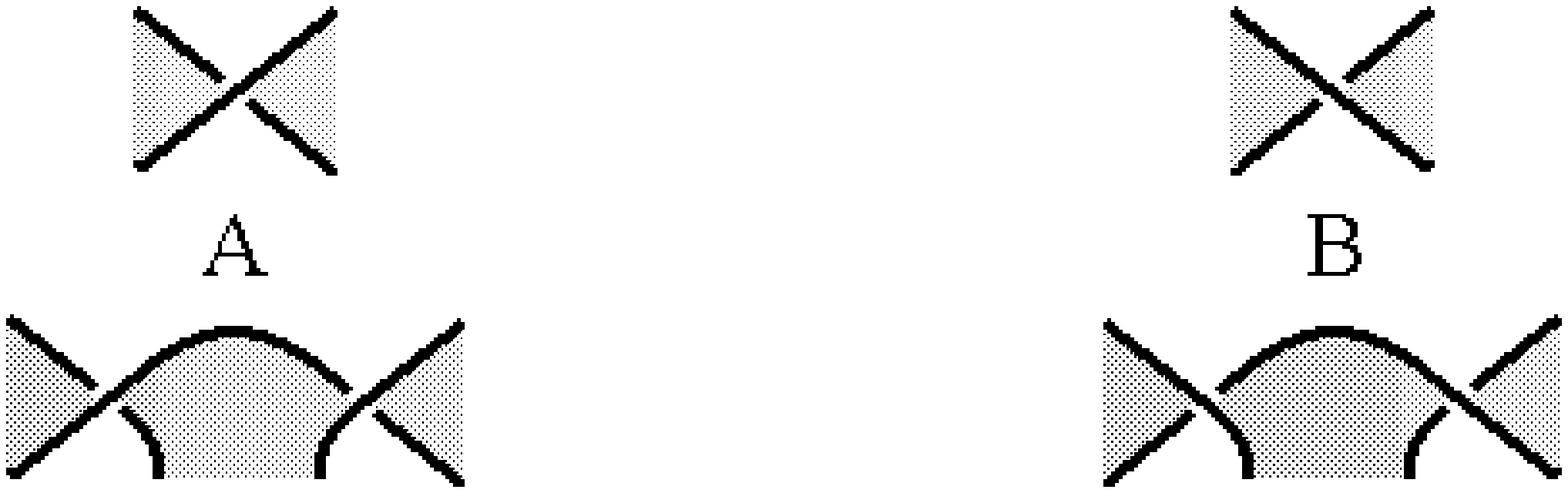} \]

As it is shown above, adjacent alternating crossings have the same coloring 
pattern of nearby regions. Hence, in a coloring of a nonsplit alternating 
diagram $D$, only one of the pattern A or B appears for every crossing. 
Reversing the coloring changes that pattern.

Resolutions of a colored diagram have induced colorings.

\[ \epsfxsize=2.20in\epsfbox{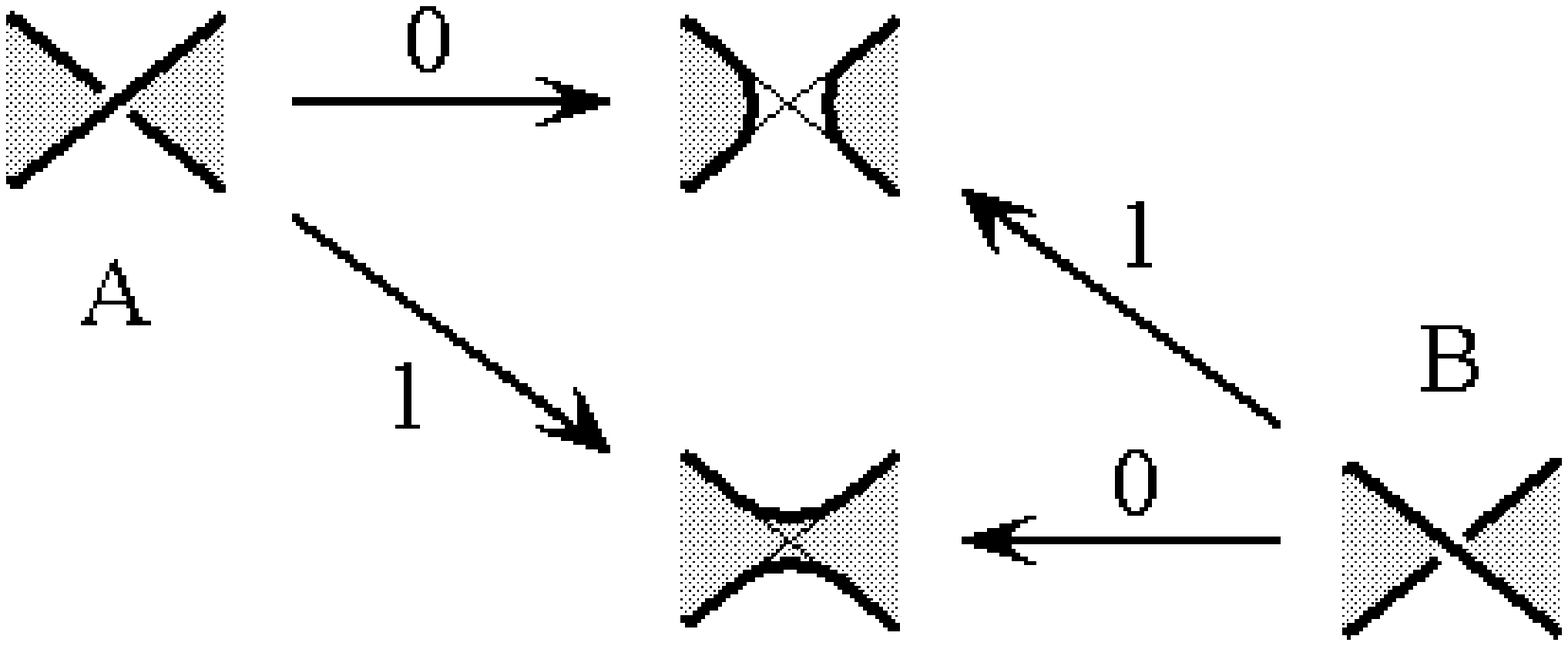} \]

\begin{defn}
For a nonsplit alternating diagram $D$, \emph{the coloring} of $D$ is the 
coloring of $D$ in which only pattern A appears. \emph{The coloring} of a 
resolution of $D$ is the coloring of that resolved diagram induced from 
the coloring of $D$.
\end{defn}

For the coloring of $D(\emptyset)$ (0-resolutions of pattern A), the trace 
of each crossing lies in a white region. Now, our claim is:

\begin{prop}
\label{prop:goodD}
For a reduced nonsplit alternating diagram $D$, the components of 
$D(\emptyset)$ bound nonoverlapping black disks in the coloring of it. 
Each black disk corresponds to each of the black regions in the coloring 
of $D$. Furthermore, every pair of black disks are connected by a chain 
of black disks, which are connected by the trace of the crossings of $D$. 
Also, no trace of crossing connects a black disk to itself.
\end{prop}

Here is visualization of our claim for diagrams of the lefthanded trefoil and 
the figure 8 knot. (The unbounded black region shown below is a disk in $S^2$.)

\[ \epsfxsize=4.03in\epsfbox{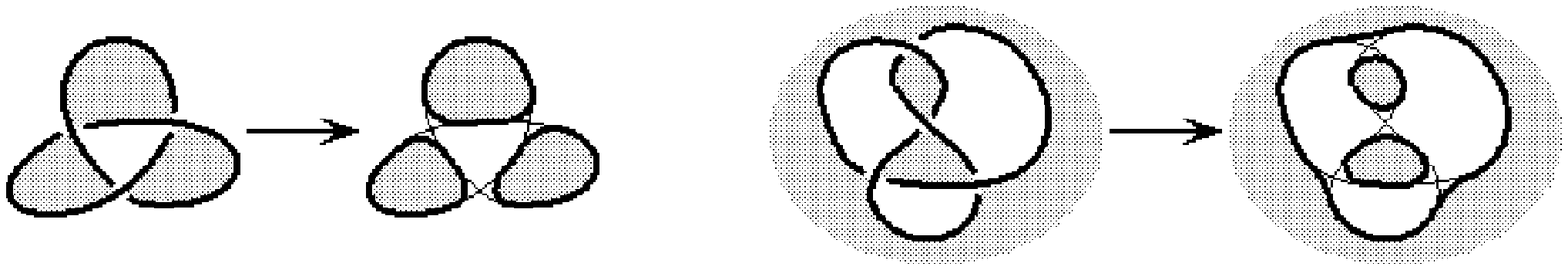} \]

\begin{proof}
At each crossing, its 0-resolution separates incident black regions. That 
gives a correspondence between the black regions in the coloring of $D$ and 
those in the coloring of $D(\emptyset)$. (While most white regions of $D$ 
merge in the process.)

In the coloring of $D(\emptyset)$, there is no trace of crossing in black 
regions. That implies:
\begin{itemize}
\item if there is a black region which is not a disk, then $D$ is split.
\item if there is a pair of black disks which cannot be connected by any 
chain, then $D$ is split.
\item if there is a trace of crossing connecting a black disk to itself, 
that crossing is removable, so $D$ is not reduced. 
\end{itemize}

\[ \epsfxsize=3.49in\epsfbox{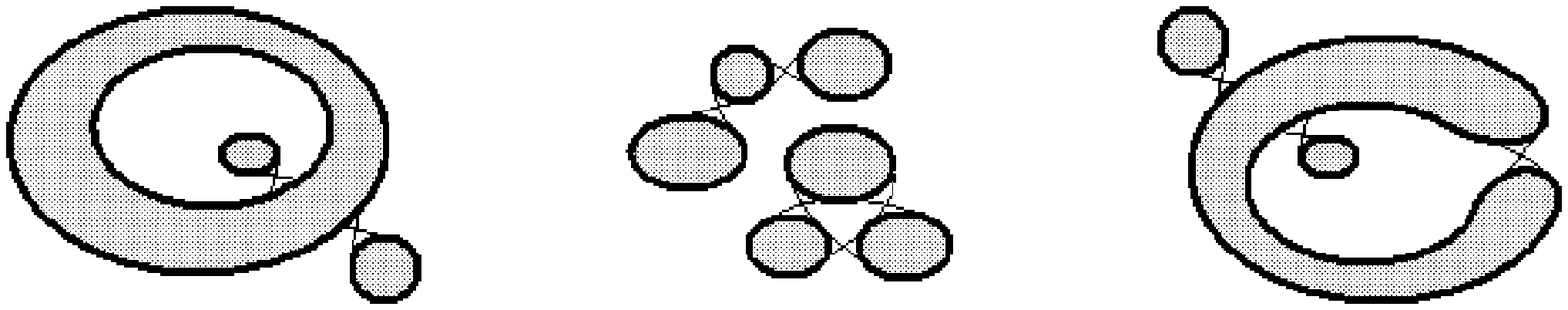} \]
\end{proof}

\begin{defn}
For a link diagram $D$, let $c(D)$ be the number of crossings of $D$, and 
$o(D)$ be the number of components of $D(\emptyset)$.
\end{defn}

For $D$ in proposition \ref{prop:goodD}, $o(D)$ agrees with the number of 
black disks in the coloring of $D(\emptyset)$.

Let $\indI$ be an ordered set of crossings of $D$. Note that $D(\indI)$ 
agrees with $D^!(\emptyset)$, and that
$ o(D) + o(D^!) $ equals the total number of black and white regions in 
the coloring of $D$, which is $ c(D) + 2 $.

We need one further step for the inductive argument to be used in our 
proof of theorem \ref{thm:conj2}.

\begin{prop}
\label{prop:betterD}
Let $D$ be a reduced nonsplit alternating link diagram with $c(D) > 0$.
Then one of the following holds.

\begin{itemize}
\item[(A)] There is a pair of black disks in the coloring of $D(\emptyset)$ 
connected by exactly one crossing.
\item[(B)] There is a pair of black disks in the coloring of $D^!(\emptyset)$ 
connected by exactly one crossing.
\item[(C)] $D$ is a connected sum of $D'$ and the Hopf link, for another 
nonsplit alternating link diagram $D'$ with $c(D) - 2$ crossings. 
\[ \epsfxsize=1.91in\epsfbox{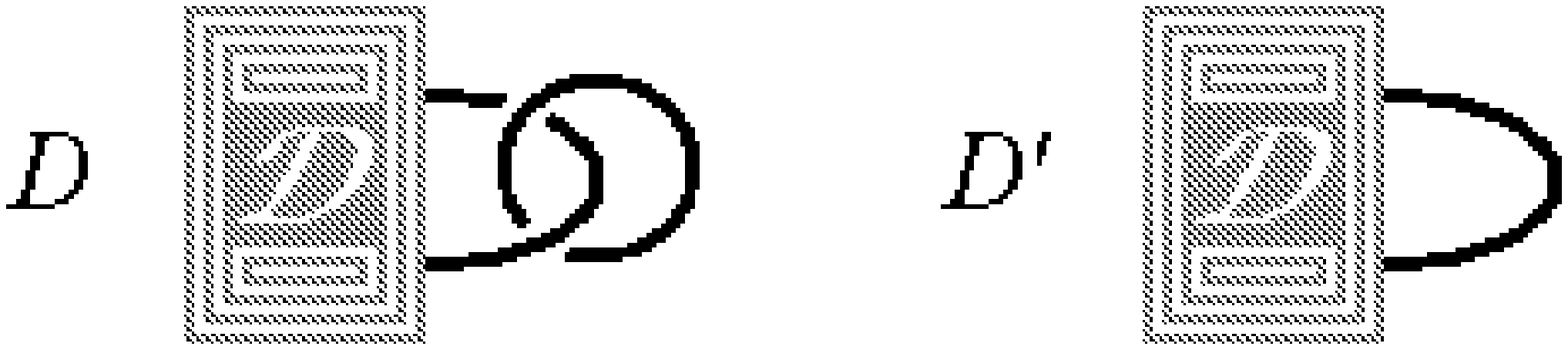} \]
\end{itemize}
\end{prop}

\begin{proof}
Since $o(D) + o(D^!) = c(D) + 2$, one of the following holds.

\begin{itemize}
\item[(a)] $o(D) > c(D)/2 + 1$.
\item[(b)] $o(D^!) > c(D)/2 + 1$.
\item[(c)] $ o(D) = o(D^!) = c(D)/2 + 1$.
\end{itemize}

[(a) $\Rightarrow$ (A)]
For $o(D)$ black disks to be connected to each other by chains of connected 
disks, there are at least $o(D) -1$ different pairs that are connected by 
crossings. If $2 (o(D) -1) > c(D) $, then at least one of those pairs is 
connected by exactly one crossing.

[(b) $\Rightarrow$ (B)]
Same as (a) $\Rightarrow$ (A).

[(c) \& not (A) \& not (B) $\Rightarrow$ (C)]
To fail (A), there are exactly $o(D) - 1$ different pairs that are connected 
by crossings and those pairs are connected by exactly two crossings.

Consider a graph consists of $o(D)$ vertices and $o(D) -1$ edges. Each vertex 
represents each black disk. For each pair of black disks connected by two 
crossings, there is an edge joining the corresponding pair of vertices.
This graph is connected, so it is a tree.

For an edge $\{a,b\}$, mark the $a$-end of it with arrow if the two crossings 
connecting the disk $a$ and $b$ are adjacent on the boundary of $a$.
For example, 

\[ \epsfxsize=3.48in\epsfbox{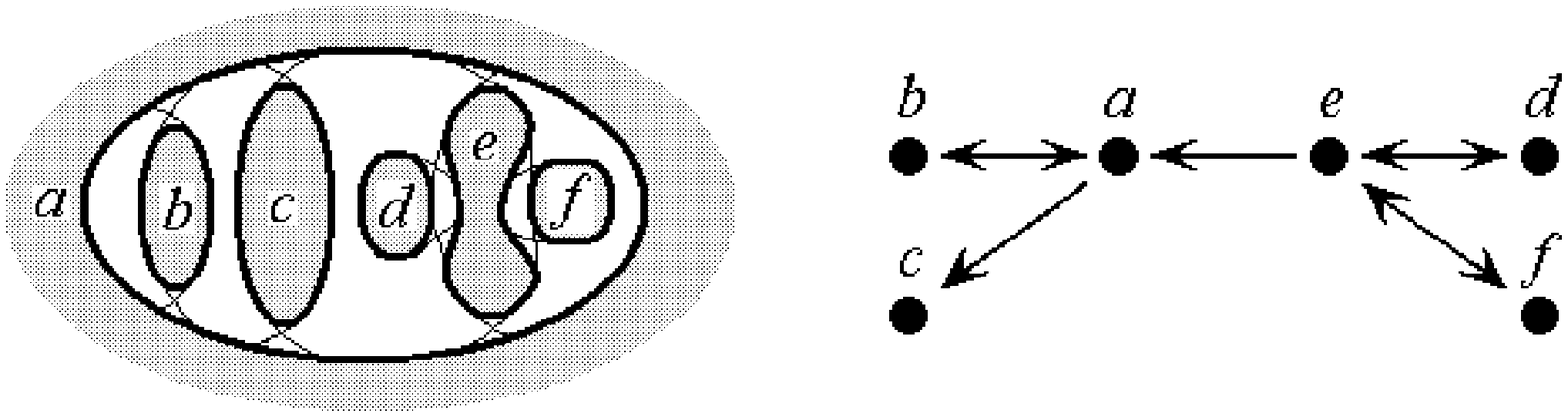} \]

A vertex of a tree is called a pendent vertex if it is incident with only one 
edge, and an edge is called a pendent edge if it is incident with a pendent 
vertex.
If $a$ is a pendent vertex, the unique edge incident with $a$ is necessarily 
marked at the $a$-end. If $b$ is not a pendent vertex, at least two edges 
have marked $b$-end, since the two crossings connecting disks $b$ and $c$ 
and those connecting disks $b$ and $d$ never alternate. 

\[ \epsfxsize=1.10in\epsfbox{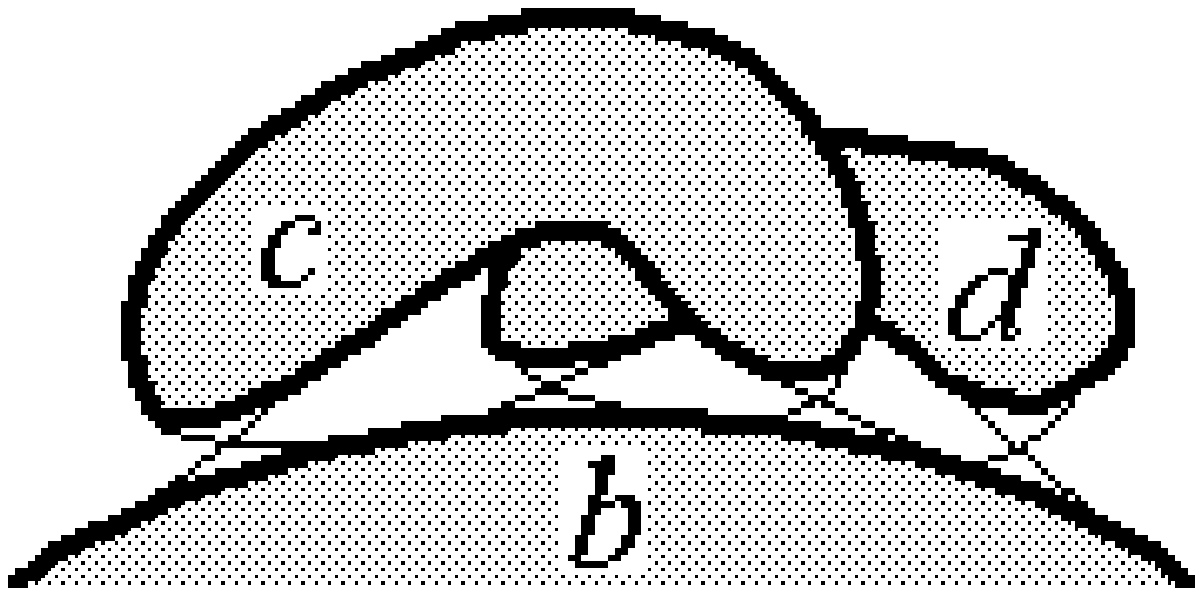} \]

If $o(D)=2$, then there is only one edge, that is a pendent edge, and both 
ends of that edge is marked. For $o(D) > 2$, let $p$ be the number of the 
pendent vertices. The number of the pendent edges is also $p$. There are 
at least $p + 2(o(D)-p)$ marked ends, but the number of nonpendent edges 
is $o(D)- 1- p$, so there is at least one pendent edge with both ends marked.
That implies (C) (up to relocation of $\infty$). 

\[ \epsfxsize=2.33in\epsfbox{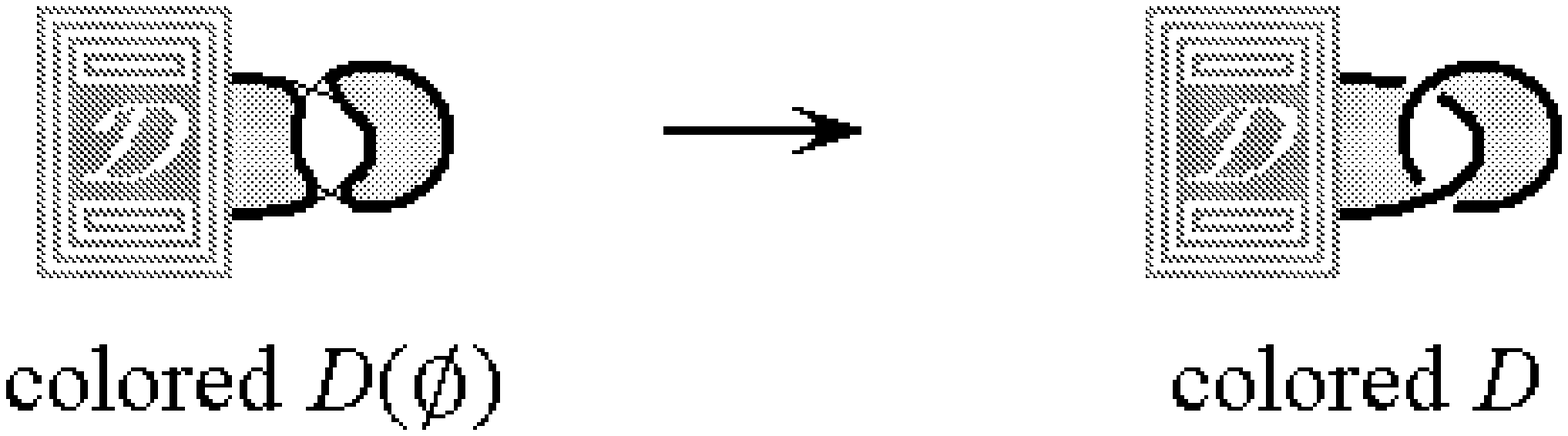} \]
\end{proof}

\begin{lemma}
\label{lemma:box-support}
For a reduced nonsplit alternating diagram $D$, $\bhm^{i,j}(D)$ is supported 
in the box $0 \leq i \leq c(D)$ and $-o(D) \leq j \leq 2c(D)-o(D)+2$,
with $\bhm^{0,-o(D)}(D) = \bhm^{c(D),2c(D)-o(D)+2}(D) =
\mathbb{Q}$.
\end{lemma}

\begin{proof}
First of all, it is clear from the construction of $\bcn(D)$ that 
$\bcn^{i,j}(D) = 0$ unless $0 \leq i \leq c(D)$.

When a resolution of $D$ is changed to another resolution of $D$ by 
replacing one 0-resolution by 1-resolution, the number of components 
either increases or decreases by one. That ensures $\bcn^{i,j}(D)$ to 
be supported in $-o(D) \leq j \leq 2c(D)-o(D)+2$.

Proposition \ref{prop:goodD} implies that $D(\emptyset)$ has one more 
component than any $D(a)$ has, because two black disks merge into one 
in the process. In terms of $\bcn^{i,j}(D)$, this means
\[ \bcn^{i,j}(D) =
\begin{cases}
\mathbb{Q} & \textrm{ if } i=0, j = -o(D) \\
0 & \textrm{ if } i > 0, j = -o(D)\textrm{ ,}
\end{cases} \]
so one half of the result follows.

For the other half, look at the other end. $D(\indI)= D^!(\emptyset)$ also 
has one more component than any $D(\indI-\{a\})= D^!(a)$ has, so that
\[ \bcn^{i,j}(D) =
\begin{cases}
\mathbb{Q} & \textrm{ if } i=c(D), j = 2c(D) -o(D)+2 \\
0 & \textrm{ if } i < c(D), j = 2c(D)-o(D)+2 \textrm{ .}
\end{cases} \]
\end{proof}

Let $D$ be a diagram satisfying (A) in proposition \ref{prop:betterD}. Let 
$a$ be a crossing of $D$ connecting a pair of black disks that no other 
crossing connects. Choose an ordering of $\indI$ in which $a$ comes the last. 
Then, $D(*0)$ still has the property that $D(*0)(\emptyset)$ has one more 
component than any $D(*0)(b)$ has. The use of (A) is that it allows $D(*1)$ 
to have that property, too.

\begin{cor}
\label{cor:box-support}
In the above setting, $\bhm^{i,j}(D(*0))$ is supported in the box 
$0 \leq i \leq c(D(*0))$ and $-o(D(*0)) \leq j \leq 2c(D(*0))-o(D(*0))+2$,
with $\bhm^{0,-o(D(*0))}(D(*0)) = \mathbb{Q}$, and
$\bhm^{i,j}(D(*1))$ is supported in the box $0 \leq i \leq
c(D(*1))$ and $-o(D(*1)) \leq j \leq 2c(D(*1))-o(D(*1))+2$, with
$\bhm^{0,-o(D(*1))}(D(*1)) =
\bhm^{c(D(*1)),2c(D(*1))-o(D(*1))+2}(D(*1)) = \mathbb{Q}$.
\end{cor}

Finally, to apply induction hypothesis to $D(*0)$ and $D(*1)$ later on, they 
need to be nonsplit alternating.

\begin{prop}
In the above setting, $D(*0)$ and $D(*1)$ are nonsplit alternating.
\end{prop}

\begin{proof}
Alternating property is easy to see.

To be nonsplit, their black disks in the induced coloring have to be 
connected. That is clear for $D(*1)$. For $D(*0)$, if the black disks 
of $D(\emptyset)$ are disconnected after removing $a$, then $a$ was a 
removable crossing in $D$, which contradicts $D$ being reduced.
\end{proof}

\subsection{Signature of an alternating link}

This section consists of the result of \cite{GL} and an application to 
alternating links, to relate the shift with the signature in theorem 
\ref{thm:conj2}.

\begin{defn}[Goeritz matrix : following \S1 of \cite{GL}]
Let $D$ be an oriented link diagram. Color the regions of $\mathbb{R}^2$ 
(or $S^2$) divided by $D$ in checkerboard fashion. Denote the white 
regions by $X_0, X_1, \cdots, X_n$. Assume that each crossing is incident 
to two distinct white regions. Assign an incidence number $\eta(a)= \pm 1$ 
to each crossing $a$ as in the figure below. For $0 \leq i , j \leq n $ define
\[ g_{ij} =
\begin{cases}
- \sum_{a \textrm{ incident to both }X_i \textrm{ and }X_j} \eta(a) 
& \textrm{ for }i \neq j \\
- \sum_{0 \leq k \leq n , k \neq i} g_{ik} & \textrm{ for }i = j \textrm{ .}
\end{cases} \]

The Goeritz matrix $G(D)$ of $D$ is the $n \times n $ 
(not $ (n+1) \times (n+1)$ !) symmetric matrix 
$G(D)=(g_{ij})_{1 \leq i ,j \leq n}$.
\end{defn}

\[ \epsfxsize=2.93in\epsfbox{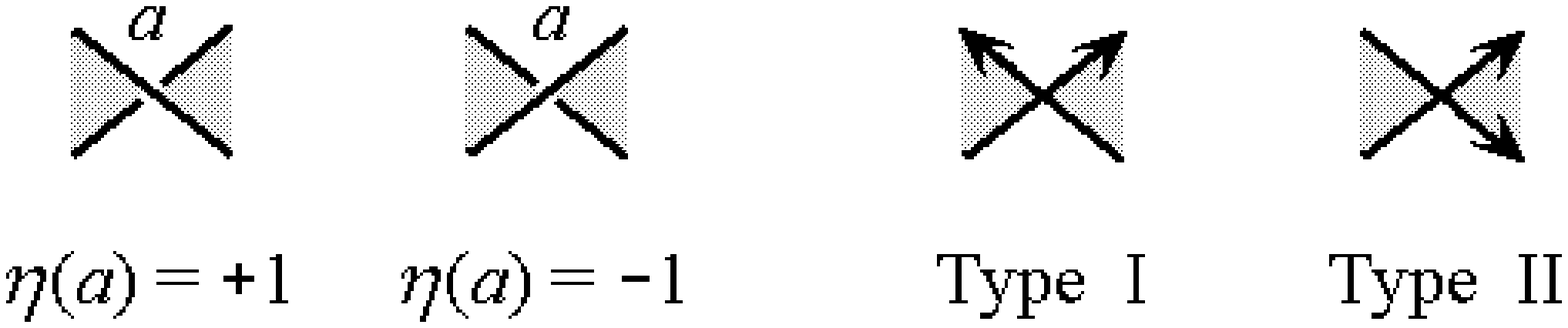} \]

The signature of an oriented link can be obtained from the signature 
of Goeritz matrix of its diagram by adding a correction term.

\begin{thm}[Theorem 6 in \cite{GL}]
For an oriented link $L$,
\[ \sigma(L) = \sgn G(D) - \mu(D) \]
for its diagram $D$, where
$ \mu(D) = \sum \eta(a) $, summed over all crossings $a$ of type II. 
(See the figure above.)
\end{thm}

\begin{prop}
\label{prop:shift}
For an oriented nonsplit alternating link $L$ and a reduced alternating 
diagram $D$ of $L$, $\sigma(L) = o(D) - y(D) -1$.
\end{prop}

\begin{proof}
In the reversed coloring of $D$, $\eta(a)= 1$ for any crossing $a$, 
components of the resolution $D(\emptyset)$ bound nonoverlapping white 
disks, $+$ crossings are of the type II, and $-$ crossings are of the 
type I, so that $g_{ij} \leq 0$ for $i \neq j$, $g_{ii} \geq 0$, and 
$\mu(D) = y(D)$. Reducedness of $D$ ensures that each crossing is incident 
to two distinct white regions.

Since
\[ \sum_{1\leq i,j\leq n}g_{ij}x_ix_j = \sum_{1 \leq i < j \leq n}
|g_{ij}|(x_i-x_j)^2 + \sum_{i=1}^n|g_{i0}|x_i^2 \geq 0, \]
$G(D)$ is a positive-definite matrix, and hence, 
\[ \sigma(L) = \sgn G(D) - \mu(D) = \rnk G(D) - y(D) = o(D) -1 - y(D) . \]
\end{proof}

\subsection{Proof of theorem \ref{thm:conj2}}

The proof is based on induction on the number of crossings of a link diagram. 
First, we prove the theorem for some number $s(L)$ instead of $\sigma(L)$, 
and then, show $s(L)=\sigma(L)$.
For convenience of proof, we will restate theorem \ref{thm:conj2} in more 
detailed and extended form as follows.

\begin{thm}
\label{thm:conj2-ext}
For any oriented nonsplit alternating link $L$, $Kh(L)(t,q)$ is supported in 
two lines $\deg(q) = 2 \deg(t) - \sigma(L) \pm 1$, its nonzero coefficient 
of the smallest degree in $t$ is on the line 
$\deg(q) = 2 \deg(t) - \sigma(L) - 1 $, its nonzero coefficient of the 
largest degree in $t$ is on the line $\deg(q) = 2 \deg(t) - \sigma(L) + 1 $, 
and those coefficients are 1.

In other words,
\[ Kh(L)(t,q) = \sum_{i=p}^m ( a_i t^i q^{2i - \sigma(L) - 1} + 
b_i t^i q^{2i - \sigma(L) + 1} ) \]
for some $p \leq m$ with $a_p = b_m = 1$.
\end{thm}

The lines $\deg(q) = 2 \deg(t) - s(L) - 1 $ and 
$\deg(q) = 2 \deg(t) - s(L) + 1 $ will be called the upper diagonal, and 
the lower diagonal, respectively, and the positions of $a_p = 1 $ and 
$b_m = 1$ will be referred to as the top at $(p, 2p - s(L) - 1)$ and 
the bottom at $(m, 2m - s(L) + 1)$, thinking of the table of coefficients 
in which the powers of $t$ increase from left to right, and the powers of 
$q$ increase from top to bottom. These terms will be applied to Khovanov's 
cohomology groups as well.

\begin{thm}
\label{thm:almost}
For any nonsplit alternating link diagram $D$, $\bhm^{i,j}(D)$ is supported 
in two lines $ j = 2i - s \pm 1 $ for some integer $s$ with the top and 
bottom on the upper diagonal and the lower diagonal, respectively.
\end{thm}

\begin{proof}
For the base case, the theorem holds for the unknotted diagram of unknot.

Assume that the statement is true for all such diagrams with less than $c$ 
crossings. Let $D$ be a nonsplit alternating link diagram with $c$ crossings. 
If $D$ is not reduced, then $\bhm(D)$ is a shift of $\bhm(D')$ for some such 
diagram $D'$ with less than $c$ crossings, so the statement is true for $D$ 
as well.

Let $D$ be reduced. By corollary \ref{cor:mirror}, it is enough to show the 
theorem for either $D$ or $D^!$. So, we may assume that $D$ has the property 
(A) or (C) in proposition \ref{prop:betterD}.

[Case (A)]
The induction hypothesis applies to $D(*0)$ and $D(*1)$. $\bhm(D(*0)) $ is 
supported in two lines with the top at $(0, -o(D(*0)))$, and $\bhm(D(*1)) $
is also supported in two lines with the top at $(0, -o(D(*1)))$.

Since $o(D)=o(D(*0))=o(D(*1))+1$, the upper diagonal and the lower diagonal 
of $\bhm(D(*0)) $ agree with those of $\bhm(D(*1))[-1]\{-1\} $.
By theorem \ref{thm:exact} and lemma \ref{lemma:box-support}, $\bhm(D) $ is 
supported in two lines with the top at $(0, -o(D))$ and the bottom at 
$(c, 2c-o(D)+2)$.

[Case (C)]
Our $D(*0)$ and $D(*1)$ are as below, and the induction hypothesis applies 
to $D'$. 

\[ \epsfxsize=3.38in\epsfbox{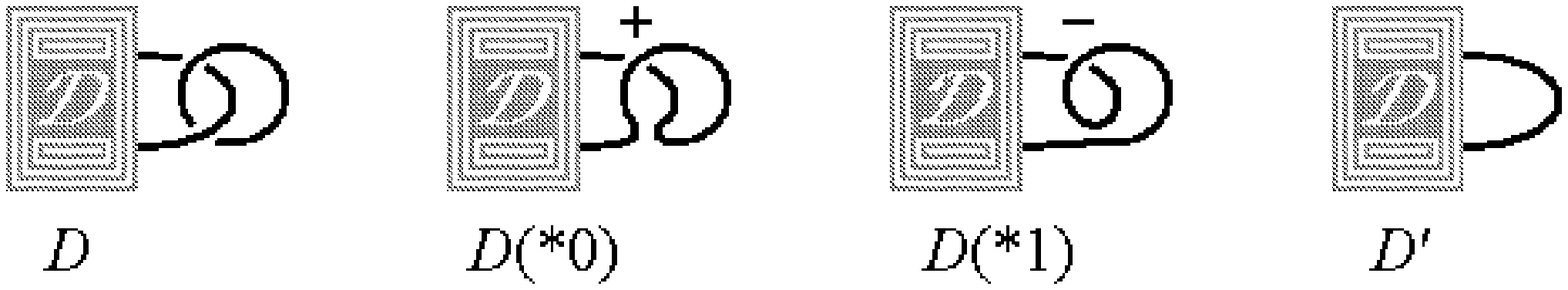} \]

Choose orientations for $D'$, $D(*0)$ and $D(*1)$ accordingly. $\bhm(D')$, 
$\bhm(D(*0))$ and $\bhm(D(*1))[-1]\{-1\}$ are shift of each other as follows.
\begin{eqnarray*}
\bhm(D(*0))
&=& \hm(D(*0))[-x(D(*0))]\{-2x(D(*0)) + y(D(*0))\} \\
&=& \hm(D')[-x(D')]\{-2x(D') + y(D')+ 1\} \\
&=& \bhm(D')[0]\{1\}
\end{eqnarray*}
\begin{eqnarray*}
\bhm(D(*1))[-1]\{-1\}
&=& \hm(D(*1))[-x(D(*1))-1]\{-2x(D(*1)) + y(D(*1)) -1\} \\
&=& \hm(D')[-x(D')-2]\{-2x(D') + y(D') -3\} \\
&=& \bhm(D')[-2]\{-3\}
\end{eqnarray*}

By induction hypothesis, $\bhm(D(*0)) $ is supported in two lines with the 
top at $(0, -o(D')-1)$, and $\bhm(D(*1))[-1]\{-1\} $
is also supported in two lines with the top at $(2, -o(D')+ 3)$. Their upper 
diagonals and lower diagonals agree.

Again, by theorem \ref{thm:exact} and lemma \ref{lemma:box-support}, 
$\bhm(D) $ is supported in two lines with the top at 
$(0, -o(D)) = (0, -o(D')-1)$ and the bottom at $(c, 2c-o(D)+2)$.
\end{proof}

Let $L$ be an oriented nonsplit alternating link and $D$ be a reduced 
alternating diagram of $L$. From theorem \ref{thm:almost}, we can conclude 
that $\hm(D) = \bhm(D)[x(D)]\{2x(D)-y(D)\} $ has the top at 
$(-x(D), -2x(D)+y(D)-o(D))$. Since the top is on the upper diagonal, our 
$s(L)$ equals $ o(D) - y(D) -1$.

In proposition \ref{prop:shift}, we saw that $\sigma(L) = o(D) - y(D) -1$. 
That finishes the proof of theorem \ref{thm:conj2}.

\bigskip

\section{An Endomorphism of the Khovanov Invariant}
\label{scn:1st}

In this section, we prove theorem \ref{thm:conj1}. The strategy of
our proof is as follows. We define a map $\my$ of degree (1,4) from
the Khovanov invariant $\hm(L)$ of any oriented link $L$ to itself,
which pairs most of $\hm(L)$. This map $\my$ added to the
coboundary map $d$ gives rise to a new cohomology theory which can
be computed explicitly. Then, we compare the cohomology groups of
$\my$ on $\hm(L)$ with the new cohomology groups of $\my + d$.

\subsection{Definition (on chain level)}
\label{sbscn:def}

Theorem \ref{thm:conj1} states that there is an almost pairing of cohomology 
groups of degree difference (1,4), so it is natural to
think of a map of degree (1,4) on the cohomology groups.

On chain level, the map $\my$ is defined in the same fashion as
the coboundary map. Instead of
\begin{eqnarray*}
\1 \otimes \1 &\buildrel m \over \longmapsto& \1 \\
\1 \otimes \x, \x \otimes \1 &\buildrel m \over \longmapsto& \x \\
\x \otimes \x &\buildrel m \over \longmapsto& 0 \\
\1 &\buildrel \Delta \over \longmapsto& \1 \otimes \x + \x \otimes \1 \\
\x &\buildrel \Delta \over \longmapsto& \x \otimes \x
\end{eqnarray*}
$\my$'s assignment is as follows.
\begin{eqnarray*}
\1 \otimes \1, \1 \otimes \x, \x \otimes \1 
&\buildrel m_{\my} \over \longmapsto& 0 \\
\x \otimes \x &\buildrel m_{\my} \over \longmapsto& \1 \\
\1 &\buildrel \Delta_{\my} \over \longmapsto& 0 \\
\x &\buildrel \Delta_{\my} \over \longmapsto& \1 \otimes \1
\end{eqnarray*}

This new multiplication $m_{\my}$ is commutative and associative.
\[ m_{\my}(m_{\my}(x \otimes y ) \otimes z) = m_{\my}(x \otimes 
m_{\my}(y \otimes z)) = 0 \textrm{ for any } x, y, z \in \A \]
The comultiplication $\Delta_{\my}$ is also cocommutative and coassociative.
\[ (\Delta_{\my} \otimes id) \circ \Delta_{\my}(z) = (id \otimes \Delta_{\my}) 
\circ \Delta_{\my}(z) = 0 \textrm{ for any } z \in \A \]
They also satisfy
\[ \Delta_{\my} \circ m_{\my} = (m_{\my} \otimes id) \circ 
(id \otimes \Delta_{\my}) \textrm{ ,} \]
since
\[ \Delta_{\my} \circ m_{\my} (y \otimes z)= (m_{\my} \otimes id) 
(y \otimes \Delta_{\my}(z)) = 0 \textrm{ for any } y, z \in \A \textrm{ .} \]
As discussed in section \ref{sbsbscn:cat}, these properties are
enough for ensuring $\my^2 = 0$.

Although $m_{\my}$ and $\Delta_{\my}$ cannot have any compatible unit or 
counit, $(\A, m_{\my}, \Delta_{\my})$ is isomorphic to 
$(\A^*, \Delta_{\my}^*, m_{\my}^*)$.

\subsection{Invariance of $\my$}
\label{sbscn:inv}

We would like to see $\my$ is well defined on $\hm(L)$. For that, $\my$ 
should (anti)commute with $d$ and be invariant under the Reidemeister moves.

\subsubsection{Anticommutativity with $d$}
\label{sbsbscn:comm}

From the viewpoint of section \ref{sbsbscn:cat}, we only need to
check the following identities.
\begin{eqnarray*}
(1) & m \circ (m_{\my} \otimes id) + m_{\my} \circ (m \otimes id) \; = \; 
m \circ (id \otimes m_{\my}) + m_{\my} \circ (id \otimes m) & \\
(2) & (\Delta \otimes id) \circ \Delta_{\my} + (\Delta_{\my} \otimes id) 
\circ \Delta \; = \; (id \otimes \Delta ) 
\circ \Delta_{\my} + (id \otimes \Delta_{\my}) \circ \Delta & \\
(3) & \Delta \circ m_{\my} + \Delta_{\my} \circ m \; = \; (m \otimes id) 
\circ (id \otimes \Delta_{\my}) + (m_{\my} \otimes id) 
\circ (id \otimes \Delta ) &
\end{eqnarray*}

\begin{proof}
(2) can be checked in the following table.

\medskip

\begin{tabular}{cccccc}
& & $(\Delta \otimes id) \circ \Delta_{\my}$ & $(\Delta_{\my} \otimes id) 
\circ \Delta$ & $(id \otimes \Delta ) \circ \Delta_{\my}$ 
& $(id \otimes \Delta_{\my}) \circ \Delta$ \\
$\1$ & $\mapsto$ & $0$ & $\1 \otimes \1 \otimes \1$ & $0$ 
& $\1 \otimes \1 \otimes \1$ \\
$\x$ & $\mapsto$ & $\1 \otimes \x \otimes \1 + \x \otimes \1 \otimes \1$ 
& $\1 \otimes \1 \otimes \x$ & $\1 \otimes \1 \otimes \x + 
\1 \otimes \x \otimes \1$ & $\x \otimes \1 \otimes \1 $
\end{tabular}

\medskip

(1) is deduced from (2) since both $(A, m, \Delta)$ and $(A, m_{\my}, 
\Delta_{\my})$ are self dual.

A table for (3) follows.

\medskip

\begin{tabular}{cccccc}
& & $\Delta \circ m_{\my}$ & $\Delta_{\my} \circ m $ & $ (m \otimes id) 
\circ (id \otimes \Delta_{\my})$ & $ (m_{\my} \otimes id) 
\circ (id \otimes \Delta )$ \\
$\1 \otimes \1$ & $\mapsto$ & $0$ & $0$& $0$& $0$\\
$\1 \otimes \x$ & $\mapsto$ & $0$ & $\1 \otimes \1$& $\1 \otimes \1$& $0$ \\
$\x \otimes \1$ & $\mapsto$ & $0$ & $\1 \otimes \1$& $0$& $\1 \otimes \1$ \\
$\x \otimes \x$ & $\mapsto$ & $\1 \otimes \x + \x \otimes \1 $ & $0$ 
& $\x \otimes \1$ & $\1 \otimes \x $
\end{tabular}

\medskip

\end{proof}

\subsubsection{Invariance under the Reidemeister moves}
\label{sbsbscn:Reid}

We also want $\my$ to commute with the isomorphisms in section
\ref{sbsbscn:move}.

[Type I]
The isomorphism was given by
\begin{eqnarray*}
\isom: &\bcn(D')& \longrightarrow \; \bcn(D)\{1\} \\
&y \otimes \1 + z \otimes \x + x& \longmapsto \; z
\end{eqnarray*}
for $y \otimes \1 + z \otimes \x \in \bcn(D'(*0)) \approx \bcn(D) \otimes \A$ 
and $x \in \bcn(D'(*1))[-1]\{-1\}$.

Then,
\begin{eqnarray*}
(\isom \my' - \my \isom)( y \otimes \1 + z \otimes \x + x)
&=& \isom (\my'_0(y) \otimes \1 + \my'_0(z) \otimes \x + 
\my'_{0 \rightarrow 1}(y \otimes \1 + z \otimes \x) - \my'_1(x)) - \my (z) \\
&=& \my(z) - \my(z) = 0\textrm{ .}
\end{eqnarray*}

[Type II]
Let $d(z)=0, z \in \bcn^i(D)[-1]\{-1\}$.
\begin{eqnarray*}
(-1)^i(\isom \my - \my' \isom) (z) &=& (-1)^i \isom \my(z) - 
\my' (z + \alpha(z)) \\
&=& -(\my'_{01}(z) + \alpha \my'_{01}(z)) - (-\my'_{01}(z) - \my'_{10} \alpha(z)
+ \my'_{01 \rightarrow 11}(z) + \my'_{10 \rightarrow 11} \alpha(z)) \\
&=& d'_{01 \rightarrow 11} \my'_{01}(z) \otimes \1 - 
\my'_{11} d'_{01 \rightarrow 11} (z) \otimes \1 - \my'_{01 \rightarrow 11}(z) \\
&=& -d'(\my'_{01 \rightarrow 11} (z) \otimes \1) \\
&& - \big( \, \bcn(D'(*11))[-2]\{-2\} \textrm{ component of } 
(d'\my'+\my'd')(z) \, \big) \otimes \1 \\
&=& -d'(\my'_{01 \rightarrow 11} (z) \otimes \1) \textrm{ .}
\end{eqnarray*}

[Type III]
If $d(x+\beta(x)+y)=0$, then $\bcn(D(*100))[-1]\{-1\}$ component of 
$d(x+\beta(x) +y)$, that is $- d_{100}(x)$, equals $0$.

In $\bcn(D)$,
\begin{eqnarray*}
\my (x +\beta(x)+ y) &=& -\my_{100}(x) +\my_{100 \rightarrow 110}(x) + 
\my_{100 \rightarrow 101}(x) - \my_{010}(d_{100 \rightarrow 110}x \otimes \1) \\
&& +\my_{010 \rightarrow 110}(d_{100 \rightarrow 110}x \otimes \1) + 
\my_{010 \rightarrow 011}(d_{100 \rightarrow 110}x \otimes \1) -\my_1 (y) \\
&=& -\my_{100}(x) +\my_{100 \rightarrow 110}(x) + \my_{100 \rightarrow 101}(x) 
- \my_{110}d_{100 \rightarrow 110}x \otimes \1 -\my_1 (y) \\
&\buildrel {(1)} \over =& -\my_{100}(x) + \beta(-\my_{100}(x)) + 
d_{110}\my_{100\rightarrow 110} x \otimes \1 + \my_{100 \rightarrow 110}(x) \\
&& + \my_{100 \rightarrow 101}(x) -\my_1 (y) \\
&\buildrel {(2)} \over \sim &
-\my_{100}(x) + \beta(-\my_{100}(x)) + \my_{100 \rightarrow 110 \approx 011} x 
+ \my_{100 \rightarrow 101}(x) -\my_1 (y) \textrm{ .}
\end{eqnarray*}
(1) is from
\begin{eqnarray*}
0&=& \big( \, \bcn(D(*110))[-2]\{-2\} \textrm{ component of } 
(d\my + \my d)(x) \, \big) \\
&=& d_{100\rightarrow 110}(-\my_{100}x) + d_{110}\my_{100\rightarrow 110} x + 
\my_{100\rightarrow 110}(-d_{100}x) + \my_{110}d_{100\rightarrow 110} x \\
\beta(-\my_{100}(x)) &=& d_{100\rightarrow 110}(-\my_{100}(x)) \otimes \1 \\
&=& - (d_{110}\my_{100\rightarrow 110} x + \my_{110}d_{100\rightarrow 110} x) 
\otimes \1 \textrm{ .}
\end{eqnarray*}
(2) is from
\begin{eqnarray*}
d \alpha(\my_{100 \rightarrow 110} x) &=&
-d_{010}(\my_{100 \rightarrow 110} x \otimes \1) + 
d_{010 \rightarrow 110}(\my_{100 \rightarrow 110} x \otimes \1) + 
d_{010 \rightarrow 011}(\my_{100 \rightarrow 110} x \otimes \1) \\
&=& -d_{110}(\my_{100 \rightarrow 110} x) \otimes \1 - \my_{100 \rightarrow 
110} x + \my_{100 \rightarrow 110 \approx 011} x \textrm{ .}
\end{eqnarray*}

Similarly, in $\bcn(D')$,
\begin{eqnarray*}
\my' ( x +\beta'(x) +y) &\sim & -\my'_{010}(x) + \beta' (-\my'_{010}(x)) - 
\my'_{010 \rightarrow 110 \approx 101}(x) + \my'_{010 \rightarrow 011}(x) - 
\my'_1 (y)\\
&=& \isom( -\my_{100}(x) + \beta(-\my_{100}(x)) + \my_{100 \rightarrow 101}(x) 
+ \my_{100 \rightarrow 110 \approx 011} x -\my_1 (y))\textrm{ .}
\end{eqnarray*}

\subsubsection{Example}

Let $T$ be the lefthanded trefoil with the diagram $D$ in section
\ref{sbsbscn:trefoil}. We have computed $\hm(T) = \hm(D)$ in
section \ref{sbsbscn:trefoil}.

Since $\my$ is of degree (1, 4), the only possible place $\my$ can be 
nontrivial is from $\hm^{-3,-9}(T)$ to $\hm^{-2,-5}(T)$. The value of 
$\my$ at a generator $[\x \otimes \x \otimes \x]$ of 
$\hm^{-3, -9}(T) = \mathbb{Q}$ is
\[ \my([\x \otimes \x \otimes \x ]) = \big[ \begin{bmatrix} \1 \otimes \x \\ 
\1 \otimes \x \\ \1 \otimes \x \end{bmatrix} \big] \textrm{ ,} \]
which is a generator of $\hm^{-2, -5}(T) = \mathbb{Q}$, so $\my$ is 
nontrivial there. 

\[ \epsfxsize=1.34in\epsfbox{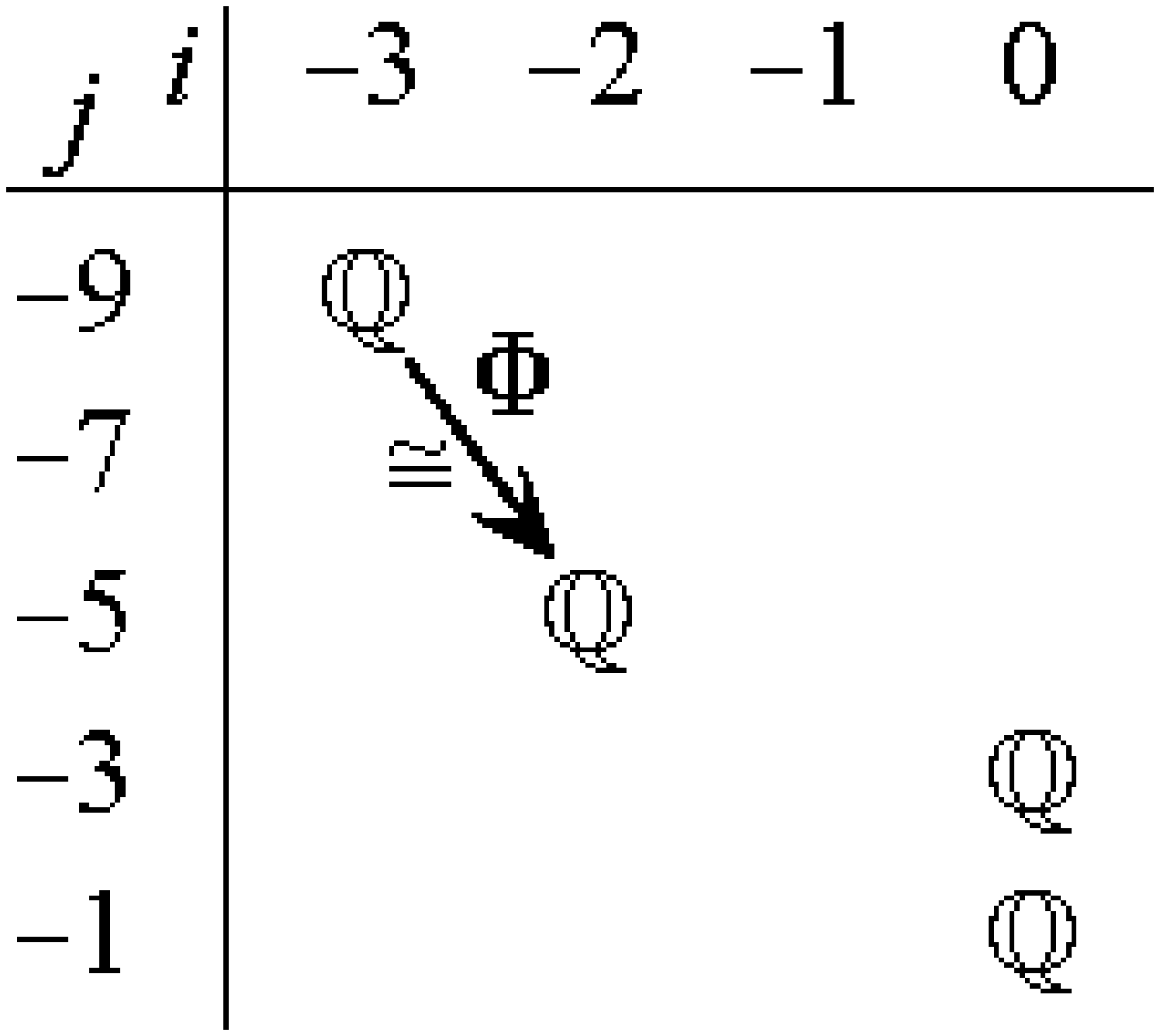} \]

\subsection{$\my + d$ and change of variables}
\label{sbscn:cv}

\subsubsection{Change of variables}

Let us forget the grading and make a change of variables as follows.
\[ \fa = \x + \1 \hspace{1in} \fb = \x - \1 \]
Since ${(\my + d)}^2 = \my^2 + \my d + d \my + d^2 = 0$, we can regard 
$\my + d$ as a new coboundary map.
\begin{eqnarray*}
\fa \otimes \fa & {\buildrel m_{\myd} \over \longmapsto} & 2 \fa \\
\fa \otimes \fb , \fb \otimes \fa & {\buildrel m_{\myd} \over \longmapsto} 
& 0 \\
\fb \otimes \fb & {\buildrel m_{\myd} \over \longmapsto} & -2 \fb \\
\fa & {\buildrel \Delta_{\myd} \over \longmapsto} & \fa \otimes \fa \\
\fb & {\buildrel \Delta_{\myd} \over \longmapsto} & \fb \otimes \fb
\end{eqnarray*}

\begin{defn}
For an oriented diagram $D$, $\nbhm (D)$ is the cohomology of the 
chain complex \linebreak 
$( \bcn(D), \my + d )$, and $\nhm (D)$ is that of $( \cn(D), \my + d )$.
\end{defn}

\subsubsection{Invariance of $\nhm (D)$ under the Reidemeister moves}
\label{sbsbscn:move2}

It would not be interesting if we can define only $\nhm (D)$, but
not $\nhm (L)$. Our proof of invariance follows chapter 5 of
\cite{K} (summarized in section \ref{sbsbscn:move}) with only the
multiplication and comultiplication maps replaced by those of $\my
+ d$. Details are left to the reader. See section
\ref{sbsbscn:move} for figures.

[Type I]
Define
\[ \tilde{X}_1 = \Ker \mydpr_{0 \rightarrow 1} \]
and
\[ \tilde{X}_2 = \{ y \otimes \frac{1}{2}(\fa-\fb) + z | y \in \bcn(D), 
z \in \bcn(D'(*1))[-1] \} \textrm{ .} \]
$\bcn(D')$ is decomposed as $\tilde{X}_1 \oplus \tilde{X}_2$ as chain 
complexes, $\tilde{X}_2$ is acyclic, and
\begin{eqnarray*}
\tilde{\isom}:& \tilde{X}_1 &\longrightarrow \; \bcn(D) \\
&y \otimes \fa + z \otimes \fb &\longmapsto \; y + z
\end{eqnarray*}
induces an isomorphism between $\nhm(D')$ and $\nhm(D)$.

[Type II]
This time,
\begin{eqnarray*}
\tilde{X}_1 &=& \{ z + \tilde{\alpha}(z) | z \in \bcn(D'(*01))[-1] \} \\
\tilde{X}_2 &=& \{ z + \mydpr y | z, y \in \bcn(D'(*00)) \} \\
\tilde{X}_3 &=& \{ z + y \otimes \frac{1}{2}(\fa-\fb) | z, 
y \in \bcn(D'(*11))[-2]\}
\end{eqnarray*}
where $\tilde{\alpha} (z) = - \mydpr_{01 \rightarrow 11}(z) \otimes 
\frac{1}{2}(\fa-\fb) \in \bcn(D'(*11))[-1] \otimes \A \approx 
\bcn(D'(*10))[-1]$.

Then, $\bcn(D') = \tilde{X}_1 \oplus \tilde{X}_2 \oplus \tilde{X}_3$ as chain 
complexes, $\tilde{X}_2$ and $\tilde{X}_3$ are acyclic, and
\begin{eqnarray*}
\tilde{\isom}: &\bcn^i(D)[-1] \approx \bcn^i(D'(*01))[-1]& \longrightarrow \; 
\tilde{X}_1 \cap \bcn^i(D') \\
&z& \longmapsto \; (-1)^i(z + \tilde{\alpha}(z))
\end{eqnarray*}
induces an isomorphism.

[Type III]
Let $\tilde{\alpha}, \tilde{\beta}, \tilde{\alpha}', \tilde{\beta}'$ be maps 
of complexes given by
\begin{eqnarray*}
\tilde{\alpha}: &\bcn(D(*110))[-2]& \longrightarrow \; \bcn(D(*010))[-1] 
\approx \bcn(D(*110))[-1] \otimes \A \\
&z& \longmapsto \; z \otimes \frac{1}{2}(\fa-\fb) \\
\tilde{\beta}: &\bcn(D(*100))[-1]& \longrightarrow \; \bcn(D(*010))[-1] \\
&z& \longmapsto \; \tilde{\alpha} \myd_{100 \rightarrow 110} (z) \\
\tilde{\alpha}': &\bcn(D'(*110))[-2]& \longrightarrow \; \bcn(D'(*100))[-1] 
\approx \bcn(D'(*110))[-1] \otimes \A \\
&z& \longmapsto \; z \otimes \frac{1}{2}(\fa-\fb) \\
\tilde{\beta}': &\bcn(D'(*010))[-1]& \longrightarrow \; \bcn(D'(*100))[-1] \\
&z& \longmapsto \; - \tilde{\alpha}' \mydpr_{010 \rightarrow 110} (z) 
\textrm{ .}
\end{eqnarray*}

$\bcn(D)$ and $\bcn(D')$ can be decomposed as below.
\begin{eqnarray*}
\bcn(D) &=& \tilde{X}_1 \oplus \tilde{X}_2 \oplus \tilde{X}_3 \\
\tilde{X}_1 &=& \{x + \tilde{\beta}(x) + y | x \in \bcn(D(*100))[-1], 
y \in \bcn(D(*1))[-1] \} \\
\tilde{X}_2 &=& \{x + \myd y | x,y \in \bcn(D(*000))\} \\
\tilde{X}_3 &=& \{\tilde{\alpha} (x) + \myd \tilde{\alpha} (y)| x, y \in 
\bcn(D(*110))[-2] \} \\
\bcn(D') &=& \tilde{Y}_1 \oplus \tilde{Y}_2 \oplus \tilde{Y}_3 \\
\tilde{Y}_1 &=& \{x + \tilde{\beta}'(x) + y | x \in \bcn(D'(*010))[-1], y \in 
\bcn(D'(*1))[-1] \} \\
\tilde{Y}_2 &=& \{x + \mydpr y | x,y \in \bcn(D'(*000))\} \\
Y_3 &=& \{\tilde{\alpha}' (x) + \mydpr \tilde{\alpha}' (y)| x, y \in 
\bcn(D'(*110))[-2] \}
\end{eqnarray*}

As before, $\tilde{X}_2, \tilde{X}_3, \tilde{Y}_2, \tilde{Y}_3$ are acyclic,
$\bcn(D(*100))[-1]$ and $\bcn(D'(*010))[-1]$, $\bcn(D(*1))[-1]$ and 
$\bcn(D'(*1))[-1]$ are naturally isomorphic, and $\tilde{X}_1$ is isomorphic 
to $\tilde{Y}_1$ via
\[ \tilde{\isom}: x + \tilde{\beta}(x) + y \longmapsto x + \tilde{\beta}'(x) 
+ y \textrm{ .} \]

\subsection{Cohomology theory of $\my + d$}
\label{sbscn:thy}

\subsubsection{Resolutions of oriented links in orientation preserving way}

Consider the resolution of an oriented link diagram in orientation preserving 
way, that is, 0-resolutions for $+$ crossings and 1-resolutions for $-$ 
crossings. This is the standard way to get a Seifert surface from a diagram 
of an oriented link. 

Consider a graph whose vertices are in 1-1 correspondence with the components 
of this resolution and whose edges connecting a pair of vertices are 
in 1-1 correspondence with the crossings connecting the corresponding pair 
of components. Since the Seifert surface obtained by the method above is 
oriented, the graph has no cycle consisting of odd number of edges, so the 
vertices of this graph can be parted into two groups in a way that the two 
end-points of each edge do not belong to the same group. Accordingly, the 
components of this resolution can be parted into two groups in a way that 
each crossing connects a component in one group to another in the other 
group. This partition does not depend on the position of $\infty$. 

\[ \epsfxsize=4.39in\epsfbox{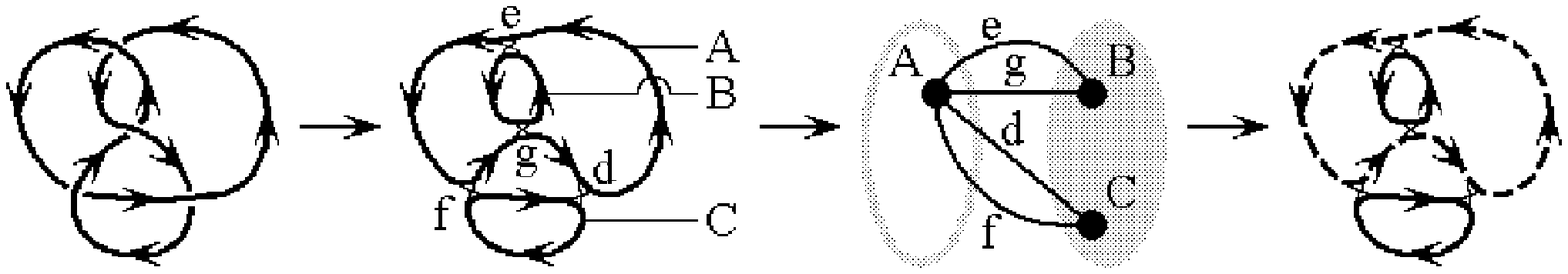} \]

\subsubsection{Hodge theory}

We can give an inner product on a chain complex so that monomials in 
$\fa, \fb$ form an orthonormal basis, then the adjoint $\myd^*$ of 
$\my + d$ is defined as follows.
\begin{eqnarray*}
\fa \otimes \fa & {\buildrel m_{\myd^*} \over \longmapsto} & \fa \\ 
\fa \otimes \fb , 
\fb \otimes \fa & {\buildrel m_{\myd^*} \over \longmapsto} & 0 \\
\fb \otimes \fb & {\buildrel m_{\myd^*} \over \longmapsto} & \fb \\
\fa & {\buildrel \Delta_{\myd^*} \over \longmapsto} & 2 \fa \otimes \fa \\
\fb & {\buildrel \Delta_{\myd^*} \over \longmapsto} & -2 \fb \otimes \fb
\end{eqnarray*}

By Hodge theory,
\[ \nbhm^i (D) \; \cong \; \Ker (\myd: \bcn^i (D) \rightarrow \bcn^{i+1} (D)) 
\cap \Ker (\myd^*: \bcn^i (D) \rightarrow \bcn^{i-1} (D))\textrm{ .} \]

\subsubsection{Computation of $\nhm$}

For an $n$ component link, there are $2^{n-1}$ different (relative) 
orientations. Each of them gives a distinct resolution when the link 
is resolved in orientation preserving way. Since no crossing connects 
components in the same group in such a resolution, the two monomials 
consisting of $\fa$ for the components in one group and $\fb$ for those 
in the other clearly belong to $\Ker \myd \cap \Ker \myd^*$. The claim 
is that these are all, i.e., others are linear combinations of these.

\begin{thm}
The dimension of $\nhm (L) = \oplus_{i \in \mathbb{Z}} \nhm^i (L)$ for 
an oriented link $L$ of $n$ components equals to $2^n$.
\end{thm}

\begin{proof}
As is in theorem \ref{thm:exact}, we have a long exact sequence of 
cohomology groups.
\[ \cdots \rightarrow \nbhm^{i-1} (D(*0)) \rightarrow \nbhm^{i-1} (D(*1)) 
\rightarrow \nbhm^i (D) \rightarrow \nbhm^i (D(*0)) 
\rightarrow \nbhm^i (D(*1)) \rightarrow \cdots \]
So, $\dim \nbhm (D)$ does not exceed $\dim \nbhm (D(*0)) + \dim \nbhm (D(*1))$.

Let us prove for knots and two component links first, using induction on the 
number of crossings.

It clearly holds for the unknot. If we have a knot with a minimal diagram $D$ 
of $c$ crossings, then one of $D(*0)$, $D(*1)$ is a knot, the other is a two 
component link, and they have one fewer crossings. Suppose $D(*0)$ is a knot. 
One of the two relative orientations of $D(*1)$ is compatible with the 
relative orientation of $D$, and the other is compatible with that of $D(*0)$. 
Then in the long exact sequence, the two generators of $\nbhm^i (D(*0))$ map 
to the two generators of $\nbhm^i (D(*1))$ coming from the relative orientation 
compatible with that of $D(*0)$. If $D(*1)$ is a knot, then the two generators 
of $\nbhm^i (D(*0))$ coming from the relative orientation compatible with 
that of $D(*1)$ map to the two generators of $\nbhm^i (D(*1))$. Hence,
\[ 2 \leq \dim \nbhm (D) \leq \dim \nbhm (D(*0)) + \dim \nbhm (D(*1)) - 4 
= 2\textrm{ .} \]

Let $D$ be a minimal diagram for a two component link with $c(D)= c$. If $D$ 
is a disjoint union of two knot diagrams $D_1$ and $D_2$, then 
$\nbhm (D) = \nbhm (D_1) \otimes \nbhm (D_2)$ with $c(D_1), c(D_2) \leq c$, 
and therefore,
\[ 4 \leq \dim \nbhm (D) \leq \dim \nbhm (D_1) \cdot \dim \nbhm (D_2) = 
4\textrm{ .} \]
If not, choose a crossing so that both $D(*0)$ and $D(*1)$ are knot diagrams. 
Then,
\[ 4 \leq \dim \nbhm (D) \leq \dim \nbhm (D(*0)) + \dim \nbhm (D(*1)) = 
4\textrm{ .} \]

An $n$ component link diagram $D$ is either a disjoint union of link diagrams 
of fewer components or can be resolved to two link diagrams of $n-1$ 
components. The proof that $\dim \nbhm(D) = 2^n$ goes the same way as above.
\end{proof}

We can tell exactly to which $\nhm^i(L)$ those generating monomials belong.

\begin{prop}
\label{prop:value-of-H}
Let $L$ be an oriented $n$ component link, $S_1, \cdots, S_n$ be its 
components, and $\ell_{jk}$ be the linking number of $S_j$ and $S_k$. Then,
\[ \dim \nhm^i(L) = 2 \cdot \big| \; \{ \; E \subset \{2,\cdots, n \} \; | 
\; (\sum_{j \in E, k \not\in E} 2 \ell_{jk}) = i \; \} \; \big| \]
\end{prop}

\begin{proof}
Let $O$ be the given (relative) orientation of $L$ with a diagram $D$, and 
$O'$ be another (relative) orientation obtained by reversing the orientations 
of $S_j, j \in E \subset \{2,\cdots, n \}$. Let $x(D), y(D)$ be the numbers 
of positive crossings and negative crossings with respect to $O$, and 
$x'(D), y'(D)$ be those with respect to $O'$.

Since a resolution in orientation preserving way is resolving $+$ crossings 
to its 0-resolutions and $-$ crossings to its 1-resolutions, those two 
monomials corresponding to $O'$ appear in $\hm^{x'(D) - x(D)}(D)$.

On the other hand, the number of negative crossings among the crossings 
between $S_j$ and $S_k$ does not change if none or both of the orientations 
of $S_j$ and $S_k$ are reversed, and if only one of them is reversed, the 
number is changed by
\begin{eqnarray*}
2 \cdot \ell_{jk} &=& \big( y - x \textrm{ among the crossings between } 
S_j \textrm{ and } S_k \big) \\
& =& \big( x' - x \textrm{ among the crossings between } S_j \textrm{ and } 
S_k \big) \textrm{ .}
\end{eqnarray*}
Therefore,
\[ x'(D) - x(D) = \sum_{j \in E, k \not\in E} 2 \ell_{jk} \textrm{ .} \]
\end{proof}

\subsection{Proof of theorem \ref{thm:conj1}}

In the previous section, we have computed $\nhm$. To prove
theorem \ref{thm:conj1}, we are going to relate $\nhm$ to
${\Ker (\my: \hm \rightarrow \hm)} / {\Img (\my: \hm
\rightarrow \hm)}$.

\begin{thm}
\label{thm:isom}
For any H-thin link $L$,
\[ \nhm (L) \cong \frac {\Ker (\my: \hm (L) \rightarrow \hm (L))}{\Img 
(\my: \hm (L) \rightarrow \hm (L))} \textrm{ .} \]
\end{thm}

\begin{proof}
For any link $L$, $(\bcn (L), d (L), \my (L))$ is a double complex up to an 
index shift. In the spectral sequence of the double complex $(\bcn (L), 
d (L), \my (L))$, the $E_2$ and $E_{\infty}$ term are isomorphic to \linebreak
${\Ker (\my: \hm (L)\rightarrow \hm (L))} / {\Img (\my: \hm (L)
\rightarrow \hm (L))}$ and $\nhm (L)$, respectively. 

If $L$ is H-thin, then $d_2$ and thereafter must be zero maps because of their 
degree. Hence,
\[ \nhm (L) \cong E_{\infty} \cong E_2 \cong \frac {\Ker (\my: \hm (L) 
\rightarrow \hm (L))}{\Img (\my: \hm (L) \rightarrow \hm (L))} \textrm{ .} \]

\end{proof}

\subsection{Extension of theorem \ref{thm:conj1} for alternating links}

Theorems \ref{thm:conj2} and \ref{thm:conj1} imply that the Khovanov 
invariant of an alternating knot determines and is determined by its Jones 
polynomial and signature. This can be extended to oriented alternating links.

Let $L$ be a link satisfying the hypothesis in proposition 
\ref{prop:value-of-H} and nonsplit alternating. We already know that two 
monomials corresponding to an orientation $O'$ belong to 
$\nhm^{\sum_{j \in E, k \not\in E} 2 \ell_{jk}}(L)$.

To find out their degrees, consider $(\fa \otimes \cdots \otimes \fa 
\otimes \fb \otimes \cdots \otimes \fb) \pm (\fb \otimes \cdots \otimes \fb 
\otimes \fa \otimes \cdots \otimes \fa)$.
\begin{eqnarray*}
(\fa \otimes \cdots \otimes \fa \otimes \fb \otimes \cdots \otimes \fb) &+& 
(\fb \otimes \cdots \otimes \fb \otimes \fa \otimes \cdots \otimes \fa) \\
&=& 2 \cdot \sum \big( \textrm{monomials in } \1 \textrm{ and } \x 
\textrm{ with even number of } \1 \big) \\
(\fa \otimes \cdots \otimes \fa \otimes \fb \otimes \cdots \otimes \fb) 
&-& (\fb \otimes \cdots \otimes \fb \otimes \fa \otimes \cdots \otimes \fa) \\
&=& (-) 2 \cdot \sum \big( \textrm{monomials in } \1 \textrm{ and } \x 
\textrm{ with odd number of } \1 \big)
\end{eqnarray*}
Degrees of monomials in $\1$ and $\x$ with even numbers of $\1$ are the same 
in (mod 4), degrees of monomials in $\1$ and $\x$ with odd numbers of $\1$ 
are also the same in (mod 4), and those two differ by 2 in (mod 4). 
Therefore, we can conclude that one of $(\fa \otimes \cdots \otimes \fa 
\otimes \fb \otimes \cdots \otimes \fb) \pm (\fb \otimes \cdots \otimes \fb 
\otimes \fa \otimes \cdots \otimes \fa)$ is mapped to upper diagonal, the 
other to lower diagonal.

Now, theorem \ref{thm:conj1} can be extended as follows.

\begin{thm}
\label{thm:conj1-ext}
For an $n$ component oriented nonsplit alternating link $L$ with 
its components $S_1, \cdots, S_n$ and linking numbers $\ell_{jk}$ of 
$S_j$ and $S_k$,
\[ Kh(L)(t,q) = q^{-\sigma (L)} \big{\{} ( q^{-1} + q) 
\big( \sum_{E \subset \{2, \cdots, n \}}( t q^2)^{\sum_{j \in E, k \not\in E} 
2 \ell_{jk}} \big) + (q^{-1} + t q^2 \cdot q) Kh'(L)(tq^2) \big{\}} \]
for some polynomial $Kh'(L)$.
\end{thm}

Hence, with the linking numbers of pairs of components provided, the Khovanov 
invariant of an oriented nonsplit alternating link determines and is 
determined by its Jones polynomial and signature.



\begin{thebibliography}{99}

\bibitem{A} L. Abrams, \emph{Two-dimensional topological quantum field 
theories and Frobenius algebras}, J. of Knot Theory and its Ramifications 
{\bf 5(5)} (1996) 569-587.

\bibitem{B} D. Bar-Natan, \emph{On Khovanov's categorification of the Jones 
polynomial}, Algebraic and Geometric Topology {\bf 2-16} (2002) 337-370.

\bibitem{B2} ---------, \emph{Khovanov homology for knots and links with 
up to 11 crossings}, preprint, 2003.

\bibitem{Br} R. A. Brualdi, \emph{Introductory combinatorics}, 3rd ed., 
Prentice-Hall, 1999.

\bibitem{G} S. Garoufalidis, \emph{A conjecture on Khovanov's invariants}, 
preprint, 2001.

\bibitem{GL} C. McA. Gordon and R. A. Litherland, \emph{On the signature of 
a link}, Invent. Math. {\bf 47} (1978) 53-69.

\bibitem{K} M. Khovanov, \emph{A categorification of the Jones polynomial}, 
Duke Math. J. {\bf 101(3)} (2000) 359-426.

\bibitem{K2} ---------, \emph{Patterns in knot cohomology I}, Experimental 
Mathematics {\bf 12(3)} (2003) 365-374.

\bibitem{L} W. B. R. Lickorish, \emph{An introduction to knot theory}, 
Graduate texts in math. {\bf 175}, Springer-Verlag, 1997.

\end{thebibliography}
\end{document}